\documentclass[11pt,a4paper,twoside]{article}
\usepackage{a4}
\usepackage{color}
\usepackage{bm, amsmath, amssymb, amsthm} 

\newcommand{\Sc}{{\mbox{\rm Sc}}}

\def\<{\langle}
\def\>{\rangle}
\def\eps{\varepsilon}
\def\NN{\mathbb{N}}
\def\ZZ{\mathbb{Z}}
\def\RR{\mathbb{R}}
\def\CC{\mathbb{C}}

\def\tr{\operatorname{Tr\,}}
\def\id{\operatorname{id\,}}
\def\Div{\operatorname{div}}

\newcommand{\ds}{\operatorname{ds}}
\newcommand{\dx}{\operatorname{dx}}
\newcommand{\dy}{\operatorname{dy}}

\newcommand{\dtau}{\operatorname{d}\!\tau}

\def\Ric{\operatorname{Ric}}
\def\vol{\operatorname{vol}}

\def\calf{\mathcal{F}}

\def\eq{\hspace*{-1.5mm}&=&\hspace*{-1.5mm}}
\def\plus{\hspace*{-1.5mm}&+&\hspace*{-1.5mm}}
\def\minus{\hspace*{-1.5mm}&-&\hspace*{-1.5mm}}
\def\dt{\partial_t}

\newtheorem{corollary}{Corollary}
\newtheorem{definition}{Definition}
\newtheorem{example}{Example}

\newtheorem{lemma}{Lemma}
\newtheorem{proposition}{Proposition}
\newtheorem{theorem}{Theorem}

\author{Vladimir Rovenski\thanks{
        E-mail: rovenski@math.haifa.ac.il.
        Supported by the Marie-Curie actions grant EU-FP7-P-2010-RG, No.~276919.}
        \ and \
        Leonid Zelenko\thanks{
        E-mail: zelenko@math.haifa.ac.il}\\
        {\small Mathematical Department, University of Haifa}
}

\title{
The mixed scalar curvature flow and harmonic foliations}
\begin{document}

\date{}

\maketitle

\begin{abstract}
We~introduce and study the~flow of metrics on a foliated Riemannian manifold $(M,g)$,
whose velocity along the orthogonal distribution
is proportional to the mixed scalar curvature, $\Sc_{\,\rm mix}$.
The~flow is used to examine the question:
When a foliation admits a metric with a given property of  $\Sc_{\,\rm mix}$
(e.g., positive or negative)\/?
We observe that the flow preserves harmonicity of foliations and yields the Burgers type equation
along the leaves for the mean curvature vector $H$ of orthogonal distribution.
If $H$ is leaf-wise conservative, then its potential obeys the non-linear heat equation
$\dt u=n\Delta_\calf\,u +(n\beta_{\mathcal D}+\Phi)\,u+\Psi^\calf_1 u^{-1}- \Psi^\calf_2 u^{-3}$
with a leaf-wise constant $\Phi$ and known functions $\beta_{\mathcal D}\ge0$ and $\Psi^\calf_i\ge0$.
 We study the asymptotic behavior of its solutions and prove that
under certain conditions (in terms of spectral parameters of leaf-wise Schr\"{o}dinger operator
$\mathcal{H}_\calf=-\Delta_\calf -\beta_{\mathcal D}\id$) there exists a unique global solution $g_t$,
whose $\Sc_{\rm mix}$ converges exponentially as $t\to\infty$ to a~leaf-wise constant.
The metrics are smooth on $M$ when all leaves are compact and have finite holonomy group.
Hence, in certain cases, there exist ${\mathcal D}$-conformal to $g$ metrics,
whose $\Sc_{\rm mix}$ is negative or positive.


\end{abstract}

\section*{Introduction}

In Introduction we discuss the question on prescribing the mixed scalar curvature of a foliation
and define the flow of leaf-wise conformal metrics depending on this kind of curvature.

\textbf{1. Geometry of foliations}.
Let $(M^{n+p},g)$ be a connected closed (i.e., compact without a~boun\-dary) Riemannian manifold,
endowed with a~$p$-dimensional foliation $\calf$,
i.e., a partition into submanifolds (called leaves) of the same dimension $p$,
and $\nabla$ the \textit{Levi-Civita connection} of $g$.
The tangent bundle to $M$ is decomposed orthogonally as $T(M)={\mathcal D}_\calf\oplus {\mathcal D}$,
where the distribution ${\mathcal D}_\calf$ is tangent to $\calf$.
Denote by $(\,\cdot\,)^\calf$ and $(\,\cdot\,)^\perp$ projections onto ${\mathcal D}_\calf$ and ${\mathcal D}$, respectively.
The second fundamental tensor and the mean curvature vector field of $\calf$~are given~by
\begin{equation*}
 h_\calf(X,Y) := (\nabla_X Y)^\perp,\qquad
 H_\calf:=\tr_{\!g} h
 \qquad(X,Y\in {\mathcal D}_\calf).
\end{equation*}
 A Riemannian manifold $(M,g)$ may admit many kinds of geometrically interesting foliations.
Totally geodesic (i.e. $h_\calf=0$) and harmonic (i.e. $H_\calf=0$) foliations
are among these kinds that enjoyed a lot of
investigation of many geometers
(see \cite{G83}, and a survey  in \cite{rov-m}).
Simple examples
are parallel circles or winding lines on a flat torus,
and a Hopf family of great circles on the 3-sphere.
 The second fundamental tensor $h$ and the integrability tensor $T$
 of the distribution ${\mathcal D}$ are defined~by
\begin{equation}\label{E-hT}
 h(X,Y) := (1/2)\,(\nabla_X Y+\nabla_Y X)^\calf,\quad
 T(X,Y):=(1/2)\,[X,\,Y]^\calf
 \quad(X,Y\in {\mathcal D}).
\end{equation}
The~mean curvature vector of ${\mathcal D}$ is given by $H=\tr_{\!g} h$.
A~foliation $\calf$ is said to be \textit{Riemannian}, or \textit{transversely harmonic},
if, respectively, $h=0,\ {\rm or}\quad H=0$.
\textit{Conformal} foliations (i.e., $h=(1/n)\,H\cdot g_{\,|{\mathcal D}}$) were introduced by Vaisman \cite{v79} as foliations admitting a transversal conformal structure.
Such foliations extend the class of Riemannian foliations.

One of the principal problems of geometry of foliations reads as follows, see \cite{rw-m}:

\textit{Given a foliation $\calf$ on a manifold $M$ and a geometric property $(P)$,
does there exist a~Riemannian metric $g$ on $M$
such that $\calf$ enjoys $(P)$ with respect to $g$?}

\noindent
Such~problems of the existence and classification of metrics on foliations (first posed explicitly by H.\,Gluck
for geodesic foliations) have been studied intensively by many geometers in the 1970's.

A~foliation is \textit{geometrically taut}, if there is a Riemannian metric making
$\calf$ harmonic.
H.\,Rumm\-ler characterized such foliations by existence of an $\calf$-closed $p$-form on $M$ that is transverse to $\calf$.
D.\,Sullivan provided a \textit{topological tautness} condition for  geometric tautness.
By the Novikov Theorem (see \cite{cc2}) and Sullivan's results, the sphere $S^3$ has no $2$-dimensional taut foliations.
In~recent decades, several tools for proving results of this sort have been developed.
Among them, one may find Sullivan's {\it foliated cycles} and new \textit{integral formulae},
see \cite{wa1} and a survey in \cite{rw-m}.

\vskip.5mm\noindent
\textbf{2. The mixed scalar curvature}.
There are three kinds of Riemannian
curvature for a~foliation: tangential, transversal, and mixed (a plane that contains a tangent vector to the foliation
and a~vector orthogonal to it is said to be mixed).
 The geometrical sense of the mixed curvature follows from the fact that for a totally geodesic
foliation, certain components of the curvature tensor, see~\cite{rov-m},
regulate the deviation of leaves along the leaf geodesics.
In general relativity, the \textit{geodesic deviation equation} is an equation involving the Riemann curvature tensor,
which measures the change in separation of neighboring geodesics or, equivalently, the tidal force experienced
by a~rigid body moving along a geodesic.
In the language of mechanics it measures the rate of relative acceleration of two particles moving forward
on neighboring geodesics.

Let $\{E_i,\,{\mathcal E}_\alpha\}_{i\le n,\,\alpha\le p}$ be a local orthonormal frame on $TM$ adapted to ${\mathcal D}$ and ${\mathcal D}_\calf$.
The \textit{mixed scalar curvature} is the following function:
 $\Sc_{\,\rm mix}(g) =\sum\nolimits_{i=1}^n\sum\nolimits_{\alpha=1}^p
 R({\mathcal E}_\alpha, E_i, {\mathcal E}_\alpha, E_i)$,
see \cite{rov-m,rw-m,wa1}.
Recall the formula,~see~\cite{wa1}:
\begin{equation}\label{eq-ran}
 \Sc_{\rm mix}(g) =\Div (H + H_\calf) +\|H\|^2 +\|H_\calf\|^2 +\|T\|^2 -\|h\,\|^2 -\|h_\calf\|^2.
\end{equation}
Integrating \eqref{eq-ran} over a closed manifold and using the Divergence Theorem,
we obtain the integral formula with the total $\Sc_{\rm mix}(g)$.
Thus, (\ref{eq-ran}) yields decomposition criteria for foliated manifolds
under constraints on the sign of $\Sc_{\,\rm mix}(g)$, see \cite{wa1} and a survey in~\cite{rov-m}.

\vskip1mm
The basic question that we address in the paper is the following:
\textit{Which foliations admit a~metric with a given property of $\,\Sc_{\,\rm mix}$ (e.g., positive or negative)}?

\begin{example}\label{Ex-03}\rm
(a) If a distribution either ${\mathcal D}$ or ${\mathcal D}_\calf$ is one-dimensional and a unit vector field $N$ is tangent to it, then
the mixed scalar curvature is simply the Ricci curvature $\Ric_{\,g}(N,N)$.
On~a~foliated surface $(M^2,g)$ this coincides with the Gaussian curvature: $\Sc_{\,\rm mix}(g)=\Ric_{\,g}(N,N)=K(g)$.

(b) For any $n\ge2$ and $p\ge1$ there exists a fibre bundle with
a closed $(n+p)$-dimensional total space and compact $p$-dimensional totally geodesic fibers,
having constant mixed scalar curvature.
 To show this, consider the Hopf fibration $\tilde\pi:S^{3}\to S^2$ of a unit sphere $S^{3}$ by great circles (closed geodesics).
 Let $(\tilde F, g_1)$ and $(\tilde B, g_2)$ be closed Riemannian manifolds
with dimensions, respectively, $p-1$ and $n-2$.
 Let $M=\tilde F\times S^{3}\times\tilde B$ has the metric product $g=g_1\times g_2$.
Then $\pi:M\to S^{2}\times\tilde B$ is a fibration with a totally geodesic fiber $\tilde F\times S^1$.
Certainly, $\Sc_{\,\rm mix}(g)\equiv2>0$.
\end{example}

\noindent
\textbf{3. Flows of metrics on foliations}.
We shall examine the basic question using evolution equations.
A~\textit{flow of metrics} on a mani\-fold is a solution $g_t$ of a differential equation $\dt g=S(g)$\,,
where the~geometric functional $S(g)$
is a~symmetric $(0,2)$-tensor usually related to some kind of curvature.
This corresponds to a dynamical system in the infinite-dimensional space of all appropriate geometric structures on a given manifold.
Denote by ${\mathcal M}(M,{\mathcal D},\,{\mathcal D}_\calf)$ the space of smooth Riemannian metrics
on $M$ such that ${\mathcal D}_\calf$ is orthogonal to ${\mathcal D}$.
Elements of ${\mathcal M}(M,{\mathcal D},\,{\mathcal D}_\calf)$ are called \textit{adapted metrics} to the pair $({\mathcal D},{\mathcal D}_\calf)$.
The notion of the ${\mathcal D}$-\emph{truncated} $(r,2)$-tensor field $S^\perp$ (where $r=0,1$)
will be helpful: $S^\perp(X_1,X_2) = S(X_1^\bot,X_2^\bot)$.
The ${\mathcal D}$-\emph{truncated metric tensor} $g^\perp$ is given~by
$g^\perp(X_1,X_2)=g(X_1,X_2)$ and $g^\perp(Y,\cdot)=0$ for all $X_i\in {\mathcal D},\ Y\in {\mathcal D}_\calf$.
For ${\mathcal D}$-\emph{conformal} adapted metrics we have $S^\perp(g)=s(g)\,g^\perp$
where $s(g)$ is a smooth function on the space of metrics on $M$.

Rovenski and Walczak \cite{rw-m} (see also \cite{rw4}) studied flows of metrics
that depend on the extrinsic geometry of codimensi\-on-one foliations, and posed the question:

\textit{Given a geometric property $(P)$, can one find an $\calf$-truncated flow $\dt g = S^\perp(g)$
on a foliation $(M,\calf)$ such that the solution metrics $g_t\ (t\ge0)$
converge to a metric for which $\mathcal{F}$ enjoys~$(P)$}?

Rovenski and Wolak~\cite{rovwol} studied ${\mathcal D}$-conformal flows of metrics
on a foliation in order to prescribe the mean curvature vector $H$ of ${\mathcal D}$.
In aim to prescribe the sign of $\Sc_{\,\rm mix}(g)$, we study the following
\textit{mixed scalar curvature flow} of metrics $g_t$, see also \cite{rz}:
\begin{equation}\label{E-GF-Ricmix-mu}
 \dt g = -2\,(\Sc_{\,\rm mix}(g) -\Phi)\,g^\perp.
\end{equation}
Here $\Phi:M\to\RR$ is a leaf-wise constant function, its value is clarified in what follows.
By Lemma~\ref{L-nablaNN} (Section~\ref{subsec:tvarb}), the flow \eqref{E-GF-Ricmix-mu}
preserves harmonic (in particular, totally geodesic) foliations.

\vskip1mm
One may ask the question:
\textit{Given a Riemannian manifold $(M,g)$ with a harmonic foliation $\calf$, when do solution metrics
$g_t$ of (\ref{E-GF-Ricmix-mu}) converge
to the limit metric $\bar g$ with $\Sc_{\,\rm mix}(\bar g)$ positive or negative}?

\begin{example}\label{Ex-surfrev}\rm

(a) Let $(M^2, g_0)$ be a surface of Gaussian curvature $K$,
endowed with a unit geodesic vector field $N$.
Certainly, (\ref{E-GF-Ricmix-mu}) reduces to the following~view:
\begin{equation}\label{E-GF-KM2}
 \dt g = -2\,(K(g)-\Phi)\,g^\perp,
\end{equation}
that looks like the normalized Ricci flow on surfaces, but uses the truncated metric
$g^\perp$ instead of~$g$.

Let $k\in C^2(M)$ be the geodesic curvature of curves orthogonal to $N$.
From \eqref{E-GF-KM2} we obtain the PDE $\dt k=K_{,x}$ (along a trajectory $\gamma(x)$ of $N$).
The above yields the \textit{Burgers equation}
\begin{equation*}
 \dt k=k_{,xx}-(k^2)_{,x},
\end{equation*}
which is the prototype for advection-diffusion processes in gas and fluid dynamics, and acoustics.
When $k$ and $K$ are known, the metrics may be recovered as $g^\perp_t=g^\perp_0\exp\,(-2\int_0^t(K(s,t)-\Phi)\ds)$.

(b) For the Hopf fibration $\pi:(S^{\,2m+1},g_{\rm can})\to\CC P^m$ of a unit sphere with fiber $S^1$,
the orthogonal distribution ${\mathcal D}$ is non-integrable while $h=0$.
By~(\ref{eq-ran}), $\Sc_{\,\rm mix}=2\,m$.
Thus, $g_{\rm can}$ on $S^{\,2m+1}$ is a fixed point in ${\mathcal M}(S^{\,2m+1},{\mathcal D},\,{\mathcal D}_\calf)$
of the flow (\ref{E-GF-Ricmix-mu}) with $\Phi=2m$.
\end{example}

\textbf{4. The nonlinear heat equation}.
The solution strategy is based on deducing from \eqref{E-GF-Ricmix-mu} the forced Burgers type equation
\begin{equation*}
 \dt H +\nabla^\calf \|H\|^2= n\nabla^\calf(\Div_\calf H) +X,
\end{equation*}
for certain vector field $X$, see Proposition~\ref{T-mainA}.
If $H$ is leaf-wise conservative, i.e., $H=-n\nabla^\calf\log\,u$ for a leaf-wise smooth function $u(x,t)>0$,
this and \eqref{eq-ran} yield the non-linear heat equation
\begin{equation}\label{Ec-dtrho1}
 (1/n)\,\dt u = \Delta_\calf\,u +(\beta_{\mathcal D}+\Phi/n)\,u+(\Psi^\calf_1/n)\,u^{-1}-(\Psi^\calf_2/n)\,u^{-3},\quad u(\,\cdot\,,0)=u_0,
\end{equation}
where functions $\beta_{\mathcal D}(x)\ge0$ and $\Psi^\calf_i(x)\ge0$ are known, and $\Delta_{\,\calf}$ is the leaf-wise Laplacian, see~\cite{cc2}.
We study the asymptotic behavior of its solutions and prove that
under certain conditions (in terms of spectral parameters of leaf-wise Schr\"{o}dinger operator $\mathcal{H}_\calf=-\Delta_\calf -\beta_{\mathcal D}\id$) the flow \eqref{E-GF-KM2} has a~unique global solution
$g_t$, whose $\Sc_{\rm mix}$ converges exponentially to a~leaf-wise constant.
The metrics are smooth on $M$ when all leaves are compact and have finite holonomy group.
Thus, in certain cases, there exist ${\mathcal D}$-conformal to $g$ metrics,
whose $\Sc_{\rm mix}$ is negative or positive.

\noindent\textbf{5. The structure of the paper}.
Section~\ref{sec:main-res} contains main results (Proposition~\ref{T-main-loc},
Theorem~\ref{T-main0} and Corollari\-es~\ref{T-main0-tg}\,--\,\ref{C-005}),
their proofs and examples for one-dimensional case and for twisted products.
These are supported by results of Section~\ref{sec:Appendix}
(Theorems~\ref{prexistest}\,--\,\ref{mainth1}) about non-linear PDE
\eqref{Ec-dtrho1} on a closed Riemannian manifold.
Throughout the paper everything (manifolds, foliations, etc.) is assumed to be
smooth (i.e., $C^{\infty}$-differentiable) and oriented.
We also assume that all the leaves of a foliation $\calf$ are compact minimal submanifolds.

\section{Main results}
\label{sec:main-res}

Based on the ``linear algebra" inequality $n\,\|h\,\|^2\ge \|H\|^2$
with the equality when ${\mathcal D}$ is totally umbilical,
we introduce the following function (a measure of ``non-umbilicity" of ${\mathcal D}$):
\begin{equation}\label{E-beta}
 \beta_{\mathcal D}:=n^{-2}\big(n\,\|h\,\|^2 -\|H\|^2\big)\ge0.
\end{equation}
For $p=1$, let $k_i$ be the principal curvatures of~${\mathcal D}$. Then
 $\beta_{\mathcal D}=n^{-2}(n\,\tau_2 -\tau_1^2)=n^{-2}\!\sum\nolimits_{i<j}(k_i-k_j)^2$.
Next lemma
allows us to reduce (\ref{E-GF-Ricmix-mu}) to
the leaf-wise PDE (with space derivatives along $\calf$ only).

\begin{lemma}
[see also \cite{rz}]\label{L-CC-riccati}
Let $\calf$ be a harmonic foliation on $(M,g)$. Then \eqref{eq-ran} reads as
\begin{equation}\label{E-RicNs}
 \Sc_{\,\rm mix}(g) = \Div_\calf H -\|H\|^2_g/n
 +\|T\|^2_g -\|h_\calf\|^2_g -n\,\beta_{\mathcal D}.
\end{equation}
\end{lemma}

\noindent\textbf{Proof}. From \eqref{eq-ran}, using $H_\calf=0$ and identity
 ${\Div}\,H={\Div}_\calf\,H -\|H\|_g^2$,
we obtain
\[
 \Sc_{\,\rm mix}(g)= {\Div}_\calf\,H -\|h\|^2_g +\|T\|^2_g-\|h_\calf\|^2_g.
\]
Substituting $\|h\|^2_g= n\,\beta_{\mathcal D}+\|H\|^2_g/n$ due to \eqref{E-beta}, we get \eqref{E-RicNs}.
\qed

\vskip1mm
We denote $\nabla^\calf f:=(\nabla f)^\calf$. Given a vector field $X$ and a function $u$ on $M$,
define the functions using the leaf-wise derivatives: the divergence
$\Div_\calf X=\sum_{\alpha=1}^{p} g(\nabla_\alpha X, {\mathcal E}_\alpha)$ and the Laplacian
$\Delta_\calf\,u=\Div_\calf(\nabla^\calf u)$.
Notice that operators $\nabla^\calf,\,\Div_\calf$ and $\Delta_\calf$ (i.e., on the leaves) are $t$-independent.

The Schr\"{o}dinger operator is central to all of quantum mechanics.
By Proposition~\ref{R-consumb} (in Section~\ref{subsec:tvarb}), the flow of metrics (\ref{E-GF-Ricmix-mu}) preserves the leaf-wise Schr\"{o}dinger operator $\mathcal{H}$, given by
\begin{equation}\label{E-Schr}
 \mathcal{H}(u)=-\Delta_\calf\,u -\beta_{\mathcal D}\,u.
\end{equation}
The spectrum of $\mathcal{H}_\calf$ on any compact leaf $F$ is an infinite sequence of isolated real eigenvalues
$\lambda^\calf_0\le\lambda^\calf_1\le\ldots\le\lambda^\calf_j\le\ldots$
counting their multiplicities, and $\lim_{\,j\to\infty}\lambda^\calf_j=\infty$.
One may fix in $L_2(F)$ an orthonormal basis of corresponding eigenfunctions $\{e_j\}$,
i.e., $\mathcal{H}_\calf(e_j)=\lambda^\calf_j e_j$.
If all leaves are compact then $\lambda^\calf_j$ are leaf-wise constant functions
and $\{e_j\}$ are leaf-wise smooth functions on~$M$.

If~the leaf $F(x)$ through $x\in M$ is compact then $\lambda^\calf_0\le0$ (since $\beta_{\mathcal D}\ge0$)
and the eigenfunction $e_0$ (called the \textit{ground state})
may be chosen positive, see Proposition~\ref{P-lambda0-one}.
 The \textit{fundamental gap} $\lambda^\calf_1-\lambda^\calf_0>0$ of $\mathcal{H}_\calf$
has mathematical and physical implications
(e.g., in refinements of Poincar\'{e} inequality and a priori estimates),
it also is used to control the rate of convergence in numerical methods of computation.
Note that the least eigenvalue of operator
$-\Delta_\calf\,u -(\beta_{\mathcal D}+\frac\Phi n)\,u$ is $\lambda^\calf_0-\frac\Phi n$.

\hskip-0.62pt
An important step in the study of evolutionary PDEs is to show short-time existence/uniqueness.

\begin{proposition}\label{T-main-loc}
Let $\calf$ be a harmonic foliation on a closed Riemannian manifold $(M, g_0)$.
Then the linearization of (\ref{E-GF-Ricmix-mu}) is the leaf-wise parabolic PDE,
hence (\ref{E-GF-Ricmix-mu}) has a unique solution $g_t$ defined on a positive time interval~$[0,t_0)$
and smooth on the leaves.
\end{proposition}

We shall say that a smooth function $f(t,x)$ on $(0,\infty)\times F$ converges to $\bar f(x)$ as $t\!\to\infty$
in $C^\infty$, if it converges in $C^k$-norm for any~$k\ge0$.
It converges {\em exponentially fast} if there exists $\omega>0$ (called the \textit{exponential rate}) such that
$\lim_{\,t\to\infty} e^{\,\omega\,t}\|f(t,\cdot)-\bar f\|_{C^k}=0$ for any $k\ge0$.

Define the domain
${U}:=\{x\in M:\ \Psi^\calf_1\Psi^\calf_2\ne0\}$
and the functions
\begin{equation}\label{E-Psi-i}
 \Psi^\calf_1:=u_0^2\,\|h_\calf\|_{g_0}^2,\quad \Psi^\calf_2:=u_0^4\,\|T\|_{g_0}^2.
\end{equation}

\begin{proposition}\label{T-mainA}
Let $\calf$ be a harmonic foliation on a Riemannian manifold $(M,g_0)$
and a family of metrics $g_t\ (0\le t<t_0)$ solve (\ref{E-GF-Ricmix-mu}). Then
\begin{equation}\label{Ec-ANtau1}
 \dt H +\nabla^\calf \|H\|^2_g= n\nabla^\calf(\Div_\calf H)
 +n\nabla^\calf(\|T\|_g^2-\|h_\calf\|_g^2 -n\,\beta_{\mathcal D}).
\end{equation}
Suppose that $H_0=-n\nabla^\calf\log u_0$ for a function $u_0>0$, then
$H_t=-n\nabla^\calf\log u$ for some positive function $u:M\times[0,t_0)$, moreover,

\noindent\
 $($i$)$ if $\Psi^\calf_2\ne0$ then $u=(\Psi^\calf_2)^{1/4}\|T\|_{g_t}^{-1/2}$,
 and the non-linear~PDE \eqref{Ec-dtrho1} is satisfied.

\noindent\
 $($ii$)$ if $\Psi^\calf_1\equiv0\equiv\Psi^\calf_2$ then the potential function $u$ may be chosen as a solution of the linear PDE
\begin{equation}\label{Ec-dtrho1b}
 (1/n)\,\dt u = \Delta_\calf\,u +\beta_{\mathcal D}\,u,\qquad u(\,\cdot\,,0)=u_0.
\end{equation}
\end{proposition}

Under certain conditions, \eqref{Ec-ANtau1} and \eqref{Ec-dtrho1b} have single-point exponential attractors.
Kirsch-Si\-mon~\cite{ks2001}
studied the forced Burgers PDE on~$\RR^n$ and proved the polynomial convergence of a solution.

Based on Proposition~\ref{T-mainA}(ii) and methods of \cite{rz}, we obtain the following.

\begin{theorem}\label{C-conf}
Let $\calf$ be a totally geodesic compact foliation with integrable orthogonal distribution
on a Riemannian manifold $(M, g_0)$, and $H_0=-n\nabla^\calf\log u_0$ for a function $u_0>0$.
Then (\ref{E-GF-Ricmix-mu}) has a~unique global solution $g_t\ (t\ge0)$ smooth on the leaves.
If $\Phi=n\,\lambda^\calf_0$ then, as $t\to\infty$, the metrics $g_t$ converge in $C^\infty$
with the exponential rate $n(\lambda^\calf_1-\lambda^\calf_0)$ to the~limit metric $\bar g$ and
\[
 \Sc_{\rm mix}(\bar g)=n\lambda^\calf_0,\qquad
 \bar H=-n\nabla^\calf\log e_0.
\]
Moreover, if the leaves have finite holonomy group, then all $g_t$ and $\bar g$ are smooth on~$M$.
\end{theorem}

\begin{example}\label{E-twisted}\rm
(a) (Example~\ref{Ex-surfrev}(a) continued). Let $({\rm T}^2, g_0)$ be a torus with Gaussian curvature $K$
and a unit vector field $N$, whose trajectories are closed geodesics.
Suppose that the curvature $k$ of ortho\-gonal
(to $N$) curves obeys $k=N(\psi_0)$ for a smooth function $\psi_0$ on ${\rm T}^2$.
In this case, $-\mathcal{H}$ in \eqref{E-Schr} coincides with the leaf-wise Laplacian,
hence $\lambda^\calf_0=0$ and $e_0={\rm const}$. We also have $\Psi=T=\beta_{\mathcal D}=0$.
By~Theorem~\ref{C-conf}, the flow of metrics (\ref{E-GF-KM2}) on ${\rm T}^2$ admits a~unique solution $g_t\ (t\ge0)$. If $\Phi=0$ then, as $t\to\infty$, the metrics converge
to a flat metric, and
$N$-curves compose a rational linear foliation.

 (b) If $\beta_{\mathcal D}=0$ then $\lambda^\calf_0=0$, see Theorem~\ref{C-conf}.
 This appear for a family $(M,g_t)=M_1\times_{f_t} M_2$ of twisted products,
i.e., the manifold $M=M_1\times M_2$ with the metrics $g_t=f_t^2 g_1 + g_2\ (t\ge0)$,
where $f_t\in C^\infty(M_1\times M_2)$ are positive functions.
The submanifolds $\{x\}\times M_2$ of a twisted product compose a totally geodesic foliation $\calf$,
while $M_1\times\{y\}$ are totally umbilical
with the leaf-wise conservative mean curvature vector $H=-n\,\nabla^\calf\log f$, see~\cite{pr}.
By \eqref{E-RicNs}, we have $\Sc_{\,\rm mix}(g) = \Div_\calf H -\|H\|^2/n$.
 If the~metrics $g_t$ solve \eqref{E-GF-Ricmix-mu}
 then $H$ obeys the Burgers type equation (see \eqref{Ec-ANtau1}  with $\beta_{\mathcal D}=0 $)
\begin{equation}\label{E-Burgers-hom}
 \dt H +\nabla^\calf \|H\|^2= n\nabla^\calf(\Div_\calf H).
\end{equation}
In this case,
$\dt f=n\,\Delta_\calf f+\Phi f$
and
the function $\tilde f:=e^{-\Phi\,t} f$ obeys the heat equation $\dt\tilde f=n\,\Delta_\calf\tilde f$.
Let $M_1$ and $M_2$ be closed and $g_t=f_t^2 g_1 + g_2$ solve \eqref{E-GF-Ricmix-mu} with $\Phi=0$,
then $g_t$ converge as $t\to\infty$ in $C^\infty$ with the exponential rate~$\lambda^\calf_1$ to the metric $\bar g=\bar f^{\,2}g_1 + g_2$, where $\bar f=\!\int_{M_2} f(0,\cdot\,,y)\dy_{\!g}$, see~\cite{rz}.
\end{example}

Define $d_{u_0,e_0}:=\min_{\,F}(u_0/e_0)/\max_{\,F}(u_0/e_0)>0$.

The central result of the work is the following.

\begin{theorem}\label{T-main0}
Let $\calf$ be a harmonic compact foliation on a closed Riemannian manifold $(M, g_0)$
and $H_0=-n\nabla^\calf\log u_0$ for a smooth function $u_0>0$. If $\Phi$ obeys the inequality
\begin{equation}\label{E-cond1}
 \Phi\ge n\lambda^\calf_0 +d_{u_0,e_0}^{\,-4}\max\nolimits_{\,F}\|T\|_{g_0}^2,
\end{equation}
then (\ref{E-GF-Ricmix-mu}) admits a unique global solution $g_t\ (t\ge0)$ smooth on any leaf $F$,
moreover, if all leaves have finite holonomy group, then $g_t$ are smooth on~$M$,
and for any $\alpha\in(0,\min\{\lambda^\calf_1{-}\lambda^\calf_0,\,2\,(\Phi/n -\lambda^\calf_0)\})$
we have the leaf-wise convergence in $C^\infty$, as $t\to\infty$, with the exponential rate $n\alpha$:
\[
 \Sc_{\,\rm mix}(g_t)\to  n\lambda^\calf_0 -\Phi\le 0,\qquad
 H_t\to-n\nabla^\calf\log e_0,\qquad h_\calf(g_t)\to0.
\]
\end{theorem}

For $T=0$, condition \eqref{E-cond1} becomes $\Phi\ge n\lambda^\calf_0$, and we have the following.

\begin{corollary}\label{T-main0-tg}
Let $\calf$ be a harmonic compact foliation on a Riemannian manifold $(M, g_0)$
with integrable normal distribution and $H_0=-n\nabla^\calf\log u_0$ for a function $u_0>0$.
If $\Phi\ge n\lambda^\calf_0$ then the claim of Theorem~\ref{T-main0} holds.
\end{corollary}

The above results are summarized (due to the basic question) in the following.

\begin{corollary}\label{C-main1}
Let $\calf$ be a harmonic compact foliation on a closed Riemannian manifold $(M, g)$
and $H=-n\nabla^\calf\log u_0$ for a smooth function $u_0>0$.

\vskip1mm
\noindent\
{\rm ($i$)} Then for any $c>d_{u_0,e_0}^{-4}\!\max\nolimits_{F}\|T\|_{g}^2$ there exists a ${\mathcal D}$-conformal to $g$ metric $\bar g$ with
 $\Sc_{\,\rm mix}(\bar g)\le -c$.

\vskip2mm
\noindent\
{\rm ($ii$)} If
 $\xi^2\|h_\calf\|^2_{g} < \xi^4\|T\|_{g}^2 +d_{u_0,e_0}^{\,-4}\max\limits_{F}\|T\|_{g}^2$,
 where $\xi={u_0}/({\tilde u_0^0 e_0})$ and $\tilde u_0^0$ is defined in Section~\ref{sec:asympttestshr},

then there exists a~${\mathcal D}$-conformal to $g$ metric $\bar g$ with $\Sc_{\,\rm mix}(\bar g)>0$.
\end{corollary}

We consider applications to a Riemannian manifold $(M,g)$ with a unit vector field~$N$
(i.e., $p=1$ or/and $n=1$).
In this case, $\Sc_{\,\rm mix}(g)$ is the Ricci curvature $\Ric_{g}(N,N)$ in the $N$-direction.

\vskip2mm
\textbf{Case $p=1$}.
Let $N$ be tangent to a geodesic foliation $\calf$,
$h$ the scalar second fundamental form and $H=\tr_{\!g} h$ the mean curvature of ${\mathcal D}=N^\perp$.
We have $h_\calf=0$, and \eqref{E-RicNs} reads
 $\Ric(N,N) = \|T\|^2 +N(H) -H^2/n -n\,\beta_{\mathcal D}$.
Let the metric evolves as, see (\ref{E-GF-Ricmix-mu}),
\begin{equation}\label{E-GF-Ricmix-2}
 \dt g=-2\,(\Ric_{\,g}(N,N)-\Phi)\,g^\perp,
\end{equation}
then $H$ obeys the PDE along $N$-curves, see \eqref{Ec-ANtau1},
 $\dt H +N(H^2)= n\,N(N(H)) +n\,\big(N(\|T\|^2) -n\,\beta_{\mathcal D}\big)$.
Suppose that $H=-n\,N(\log u_0)$ for a leaf-wise smooth function $u_0>0$ on $M$,
then we assume $H=-n\,N(\log u)$
for a~positive function $u:M\times[0,t_0)$, see Proposition~\ref{T-mainA}.

If ${\mathcal D}$ is integrable, then the function $u(\,\cdot\,,\,t)>0$ may be chosen as a~solution
of the following linear heat equation, see~\eqref{Ec-dtrho1b},
 $\dt u = n\,N(N(u)) +n\,\beta_{\mathcal D}\,u$, where $u(\,\cdot\,,0)=u_0$.
By Theorem~\ref{T-main0}, the flow (\ref{E-GF-Ricmix-2}) admits a~unique global solution $g_t\ (t\ge0)$.
If $\lambda^\calf_0-\Phi/n<0$ then we have exponential convergence as $t\to\infty$ of $g_t\to \bar g$,
$H\to -n\,N(\log e_0)$ and $\Ric_{\,g_t}(N,N)\to n\lambda^\calf_0-\Phi$.

If ${\mathcal D}$ is nowhere integrable, then $u=(\Psi^\calf_2)^{1/4}\|T\|^{-1/2}_{g_t}$
(with $\Psi^\calf_2:=u_0^4\,\|T\|^2_{g_0}>0$),
moreover, the potential function $u>0$ solves the non-linear heat equation, see~\eqref{Ec-dtrho1},
\begin{equation*}
 (1/n)\,\dt u = N(N(u)) +(\beta_{\mathcal D}+\Phi/n)\,u -(\Psi^\calf_2/n)\,u^{-3},\qquad u(\,\cdot\,,0)=u_0.
\end{equation*}
If \eqref{E-cond1} are satisfied,
then (\ref{E-GF-Ricmix-2}) admits a~unique solution $g_t\ (t\ge0)$.
We have exponential convergence as $t\to\infty$ of functions
$H\to -n\,N(\log e_0)$ and $\Ric_{\,g_t}(N,N)\to n\lambda^\calf_0-\Phi$.
By Theorem~\ref{T-main0} we have

\begin{corollary}\label{C-004}
Let $N$ be a~unit vector field tangent to a geodesic foliation $\calf$ on $(M,g)$.

\noindent\ {\rm ($i$)}
Then for any $c>d_{u_0,e_0}^{\,-4}\max\nolimits_{\,F}\|T\|_{g}^2$ there is ${\mathcal D}$-conformal to $g$
metric $\bar g$ with the property
 $\Ric_{\,\bar g}(N,N)\le -c<0$.

\noindent\  {\rm ($ii$)} If~$\,\xi^2\|h_\calf\|^2_{g} <\xi^4\|T\|_{g}^2
 +d_{u_0,e_0}^{\,-4}\max_{\,F}\|T\|_{g}^2$,
 where $\xi={u_0}/({\tilde u_0^0 e_0})$  and $\tilde u_0^0$ is defined in Section~\ref{sec:asympttestshr},
then there exists ${\mathcal D}$-conformal to $g$ metric $\bar g$ such that $\Ric_{\,\bar g}(N,N)>0$.
\end{corollary}

\vskip1mm
\textbf{Case $n=1$}.
Let
$N$ be orthogonal to a compact harmonic foliation $\calf$ of codimension one.
Then $T=\beta_{\mathcal D}=H_\calf=0$, $H=\nabla_NN$,
$\Psi^\calf_2=0$,
$\Psi^\calf_1=u_0^2\|h_\calf\|_{g_0}^2$,
the operator \eqref{E-Schr} coincides with $-\Delta_\calf$
(hence, $\lambda^\calf_0=0$ and $e_0={\rm const}$),
and \eqref{E-RicNs} reads
 $\Ric(N,N) = \Div_\calf H -\|H\|^2 +\|h_\calf\|^2$.
By \eqref{Ec-ANtau1} we~have
\begin{equation*}
 \dt H +\nabla^\calf \|H\|_{g_t}^2= \nabla^\calf(\Div_\calf H)
 -\nabla^\calf(\|h_\calf\|_{g_t}^2).
\end{equation*}
Suppose the condition $H_0=-\nabla^\calf\log u_0$ for a leaf-wise smooth function $u_0>0$ on $M$.
Then $H=-\nabla^\calf\log u$, where, see \eqref{Ec-dtrho1},
\begin{equation*}
 \dt u = \Delta_\calf\,u
 +\Phi\,u+\Psi^\calf_1\,u^{-1},\qquad u(\,\cdot\,,0)=u_0.
\end{equation*}
If $\Phi>0$
then \eqref{E-cond1} holds and, by Theorem~\ref{T-main0}, the flow (\ref{E-GF-Ricmix-2})
admits a~unique global solution $g_t\ (t\ge0)$.
As $t\to\infty$, we have convergence
 $H\to 0,\ \Ric_{\,g_t}(N,N)\to -\Phi,\ h_\calf(g_t)\to0$
with the exponential rate $\alpha$ for any
$\alpha\in(0,\,\min\{\lambda^\calf_1,\,2\,\Phi\})$.
By Theorem~\ref{T-main0}, we have the following.

\begin{corollary}\label{C-005}
Let $\calf$ be a codimension one harmonic foliation with a unit normal vector field $N$.
Then for any $c>d_{u_0,e_0}^{\,-4}\max\nolimits_{\,F}\|T\|_{g}^2$ there is ${\mathcal D}$-conformal to $g$
metric $\bar g$ with
 $\Ric_{\,\bar g}(N,N)\le -c$.
\end{corollary}

For a totally geodesic foliation $\calf$, i.e., $h_\calf\equiv0$, \eqref{E-RicNs} reads
 $\Ric(N,N) = \Div_\calf H -\|H\|^2 = \Div H$.

Let the metric evolves~by \eqref{E-GF-Ricmix-2}.
By Proposition~\ref{T-mainA}, $H$ obeys the homogeneous Burgers equation
$\dt H +\nabla^\calf \|H\|^2= \nabla^\calf(\Div_\calf H)$.
Suppose that the curvature vector $H$ of $N$-curves is leaf-wise conservative:
$H=-\nabla^\calf\log u$
for a function $u>0$.
This yields the heat equation  $\dt u = \Delta_\calf\,u$.
Solution of above PDE
satisfies on the leaves
 $\bar u:=\lim\limits_{t\to\infty}u(t,x)=\int_{F_x} u_0(x)\dx/{\rm Vol}(F_x,g)$.
Since $\nabla^\calf e_0=0$, we have $\bar H=\lim\nolimits_{\,t\to\infty} H(t,\,\cdot)=0$.
Then $\Ric_{\bar g}(N,N)=0$, where $\bar g=\lim\nolimits_{\,t\to\infty} g_t$.

\vskip1mm
\textbf{Surfaces: $n=p=1$}.
 Let $(M^2, g)$ be a surface with a geodesic unit vector field $N$.
The metric in biregular foliated coordinates $(x,\theta)$ is $g=d x^2+\rho^2\,d\theta^2$,
where $\rho$ is a positive function and $\partial_x$ is the $N$-derivative.
Let $M^2\subset\RR^3: r=[\rho(x)\cos\theta,\ \rho(x)\sin\theta,\ h(x)]$ be a rotational surface,
where $0\le x\le l$, $|\theta|\le\pi,\ \rho\ge0$ and $(\rho^\prime)^2+(h^\prime)^2=1$.
 Its metric belongs to warped products, see~Example~\ref{E-twisted}.
The profile curves $\theta={\rm const}$ are geodesics tangent to $N$.

Let $M_t^2\subset\RR^3$ be a one-parameter family of surfaces of revolution
(foliated by profile curves) such that the induced metric $g_t$ obeys~(\ref{E-GF-KM2}).
The profile of $M^2_0$ (parameterized as above)
is $XZ$-plane curve $\gamma_0=r(\cdot,0)$, and $\theta$-curves are circles in $\RR^3$.
Thus $N=r_{,\,x}$ is the unit normal to $\theta$-curves on $M_t^2$.
Since the geodesic curvature of parallels is $k=-(\log\rho)_{,x}$, we have
 $\dt\rho=\rho_{,xx}-\Phi\,\rho$.

 When $\Phi=0$, the flow of metrics $\dt g =-2 K\,g^\perp$ reduces to $\rho_{,\,t} = -K\rho$.
Since $(\rho_{,\,x})_{,\,t}=(\rho_{,\,x})_{,\,xx}$,
by the maximum principle, we have the inequality $|\rho_{,\,x}|<1$ for all $t\ge0$.
When such a solution $\rho(x,t)\ (t\ge0)$ is known, we find
 $h = \int\sqrt{1-(\rho_{,\,x})^2}\dx$.
Suppose that the boundary conditions are
 $\rho(0,t)=\rho_1$, $\rho(l,t)=\rho_2$ and $h(0,t)=h_1$,
where $\rho_2\ge\rho_1\ge0$ and $t\ge0$. By the heat equation theory,
the solution $\rho$ approaches as $t\to\infty$ to a linear function
$\bar\rho=x\rho_1+(l-x)\rho_2>0$.
Also, $h$ approaches as $t\to\infty$ to a linear function $\bar h=x h_1+(l-x)h_2$,
where $h_2$ may be determined from the equality $(\rho_2-\rho_1)^2+(h_2-h_1)^2=l^2$.
The curves $\gamma_t$ are isometric one to another for all $t$ (with the same arc-length parameter $x$).
The limit curve $\lim\nolimits_{\,t\to\infty}\gamma_t=\bar\gamma=[\bar\rho,\,\bar h]$
is a line segment of length $l$. Thus, $M_t$ approach as $t\to\infty$ to the
flat surface of revolution $\bar M$ -- the patch of a cone or a cylinder generated by $\bar\gamma$.

\section{Proof of main results}
\label{sec:egf}

\subsection{Holonomy of a compact foliation}
\label{subsec:prel}


The notion of \textit{holonomy} uses that of a germ of a locally defined diffeomorphism (i.e., an equivalence class of certain maps).
The germs of diffeomorphisms $(\RR^n,0)\to(\RR^n,0)$ with fixed origin form a group, denoted by ${\rm Diff}_0(\RR^n)$.
We denote by ${\rm Diff}^+_0(\RR^n)$ the subgroup of germs of diffeomorphisms which preserve orientation of $\RR^n$.
Let $(M,\calf)$ be a foliated manifold, $x,y$ be two points on a leaf $F$ and $Q_x, Q_y$ be transversal sections (diffeomorphic to $\RR^n$).
To any path $\alpha$ from $x$ to $y$ in $F$ we associate
a germ of a diffeomorphism $h(\alpha):(Q_x,x)\to(Q_y,y)$ called the holonomy of the path,
for $x=y$, this is a generalization of
a first return map.
We obtain a homomorphism of groups $h: \pi_1(F)\to{\rm Diff}_0(\RR^n)$.
The~image ${\rm Hol}(F):=h(\pi_1(F))$ is called the \textit{holonomy group} of~$F$.
 If $\calf$ is transversally orientable then
 ${\rm Hol}(F)$ is a subgroup of ${\rm Diff}^+_0(\RR^n)$.

Certainly, ${\rm Hol}(F)$ is finite when the first fundamental group, $\pi_1(F)$, is finite.
Note that a~foliation whose leaves are the fibers of a fibre bundle has trivial holonomy group.

\begin{definition}\rm
A \textit{foliated bundle} is a fibre bundle $p:E\to F$
admitting a foliation $\calf$ whose leaves meet transversely all the fibers
$E_x=p^{-1}(x)\ (x\in F)$, and the bundle projection restricted to each leaf $F'\in\calf$ is a covering map $p:F'\to F$.
There is a representation $h_x: \pi_1(F,x)\to{\rm Diff}(E_x)$
of the first fundamental group $\pi_1(F,x)$ in the group of diffeomorphims of a fiber, called
the \textit{total holonomy homomorphism} for the foliated bundle.
\end{definition}

\textbf{Local Reeb stability Theorem} (see \cite[Vol.~I]{cc2}).
\textit{Let $F$ be a compact leaf of a foliated manifold $(M,\calf)$ and \,${\rm Hol}(F)$ is finite.
Then there is a normal neighborhood ${\rm pr} :V\to F$ of $F$ in $M$ such that $(V, \calf_{\,|V}, {\rm pr})$ is
a foliated bundle with all leaves compact (and transversal to the fibers).
Furthermore, each leaf $F'\subset V$ has
finite holonomy group of order at most the order of ${\rm Hol}(F)$ and the covering
${\rm pr}_{\,|F'}: F'\to F$ has $k$ sheets, where $k \le {order\ of}\ {\rm Hol}(F)$}.

\vskip1mm
In other words, for a compact leaf $F$ with finite holonomy group, there exists
a saturated neighborhood $V$ of $F$ in $M$ and
a diffeomorphism from $E:=\tilde F\times_{{\rm Hol}(F)}\RR^n$
under which $\calf_{\,|V}$ corresponds to the bundle foliation on $E$.
Here $\tilde F$ is a covering space of $F$ associated with ${\rm Hol}(F)$.

The following method construction foliations is related to local Reeb stability Theorem.
Suppose that a group $G$ acts freely and properly discontinuously on a connected manifold $\tilde F$
such that $\tilde F/G=F$. Suppose also that $G$ acts on a manifold $Q$.
Now form the quotient space $E:=\tilde F\times_G Q$, obtained from the product space
$\tilde E:=\tilde F\times Q$ by identifying $(gy, z)$ with $(y, gz)$ for any $y\in \tilde F,\, g\in G$ and $z\in Q$.
Thus $E$ is the orbit space of $\tilde E$ w.r.t. a properly discontinuous action of $G$.
It is also Hausdorff, so it is a manifold.
 The projection ${\rm pr}_1: \tilde E\to\tilde F$ induces a submersion ${\rm pr}:E\to F$, so we have the commutative diagram $(\tilde F\to F)\circ {\rm pr}_1={\rm pr}\circ(\tilde E\to E)$.
The map ${\rm pr}$ has the structure of a~fibre bundle over $F$ with vertical fiber $Q$.
(Fibre bundles which can be obtained in this way are exactly those with discrete structure group).
We claim that $E$ admits also horizontal leaves, so that ${\rm pr}$ maps each leaf to $F$ as a covering projection.
 Indeed, the foliation $\tilde\calf$ on $\tilde E$, which is given by the submersion
${\rm pr}_2:\tilde E\to Q$, is invariant under the action of $G$,
and hence we obtain the quotient foliation $\calf=\tilde\calf/G$ on $E$.
If $z\in Q$ and $G_z\subset G$ is the isotropy group at $z$ of the action by $G$ on $Q$,
then the leaf of $E$ obtained from the leaf $\tilde F\times\{z\}$ is naturally diffeomorphic to
$\tilde F/G_z$, and ${\rm pr}$ restricted to this leaf is the covering $\tilde F/G_z\to F$.

\subsection{${\mathcal D}$-conformal adapted variations of metrics}
\label{subsec:tvarb}

The \textit{Levi-Civita connection} $\nabla$ of a metric $g$ on $M$ is given by well-known formula
\begin{eqnarray}\label{eqlevicivita}
\nonumber
 2\,g(\nabla_X Y, Z) \eq X(g(Y,Z)) + Y(g(X,Z)) - Z(g(X,Y)) \\
 \plus g([X, Y], Z) - g([X, Z], Y) - g([Y, Z], X)\qquad(X,Y,Z\in TM).
\end{eqnarray}
Let $g_t$ be a smooth family of metrics on $(M,\calf)$ and $S=\dt g$.
Since the difference of two connections is a tensor, $\dt\nabla^t$ is a $(1,2)$-tensor on $(M,g_t)$.
Differentiating \eqref{eqlevicivita} with respect to $t$ yields, see~\cite{rw-m},
\begin{equation}\label{eq2}
 2\,g_t((\dt\nabla^t)(X, Y), Z)=(\nabla^t_X S)(Y,Z)+(\nabla^t_Y S)(X,Z)-(\nabla^t_Z S)(X,Y)
\end{equation}
for all $t$-independent vector fields $X,Y,Z\in\Gamma(TM)$.
If $S=s(g) g^\perp$, for short we write
\begin{equation}\label{E-sgeneral}
 \dt g=s\,g^\perp
\end{equation}
for a certain $t$-dependent function $s$ on $M$.
In this case, the volume form $d\vol$ evolves as \cite{ck1}
\begin{equation}\label{E-dtdvol}
 ({d}/{dt})(d\vol_t)=(n/2)\,s_t\,d\vol_t.
\end{equation}

\begin{lemma}\label{L-nablaNN}
For ${\mathcal D}$-conformal adapted variations \eqref{E-sgeneral}  of metrics we have
\begin{eqnarray}\label{Es-S-F}
 \dt h_\calf \eq -s\,h_\calf,\qquad
 \dt H_\calf = -s\,H_\calf.
\end{eqnarray}
Hence, variations \eqref{E-sgeneral} of metrics preserve harmonic and totally geodesic foliations.
\end{lemma}

\noindent\textbf{Proof}.
Let $g_t\ (t\ge0)$ be a family of metrics on $(M,\calf)$ such that $\dt g_t=S(g)$,
where the tensor $S(g)$ is ${\mathcal D}$-truncated.
Using (\ref{eq2}), we find for $X\in {\mathcal D}$ and $\xi,\eta\in {\mathcal D}_\calf$,
\begin{eqnarray*}
  2\,g_t(\dt{h_\calf}(\xi,\eta),X)\eq g_t(\dt(\nabla^t_{\xi}\,\eta)+\dt(\nabla^t_{\eta}\,\xi),\ X)
  =(\nabla^t_{\xi} S)(X,\eta)+(\nabla^t_{\eta} S)(X,\xi)-(\nabla^t_X S)(\xi,\eta)\\
 \eq-S(\nabla^t_{\xi}\,\eta, X) -S(\nabla^t_{\eta}\,\xi, X)
        =-2\,S({h_\calf}(\xi,\eta), X).
\end{eqnarray*}
Assuming $S(g)=s(g)\,g^\perp$, we have \eqref{Es-S-F}$_1$. Tracing this we have \eqref{Es-S-F}$_2$.
By the theory of ODEs, if $H_\calf=0$ or $h_\calf=0$ at $t=0$ then respectively $H_\calf=0$ or $h_\calf=0$ for all $t>0$.
\qed

\vskip1mm
The~\textit{co-nullity operator} is defined~by
 $C_N(X)=-(\nabla_{\!X} N)^\bot$, for $X\in {\mathcal D},\ N\in {\mathcal D}_\calf$.
One may decompose $C$ into symmetric and skew-symmetric parts as
 $C_N=A_N+T_N^\sharp$.
The~Weingarten operator $A_N$ of ${\mathcal D}$
and the operator $T^\sharp_N$ are related with tensors $h$ and $T$, see \eqref{E-hT},~by
\begin{equation*}
 g(A_N(X),\,Y) = g(h(X,Y),\ N),\quad
 g(T^\sharp_N(X),\,Y) = g(T(X,Y),\ N),\quad X,Y\in {\mathcal D}.
\end{equation*}
The proof of the next lemma is based on (\ref{eq2}) with $S=s\,g^\perp$.

\begin{lemma}[see \cite{rz} and \cite{rovwol}]\label{L-btAt2}
For ${\mathcal D}$-conformal adapted variations \eqref{E-sgeneral}  of metrics we have
\begin{equation}\label{Es-S-b}
 \dt A_N = -\frac12\,N(s)\,\widehat\id,\ \
 \dt T^\sharp_N = -s\,T^\sharp_N\ (N\in{\mathcal D}_\calf),\ \
 \dt H = -\frac n2\,\nabla^\calf s,\ \
 \dt(\Div_\calf H) = -\frac n2\,\Delta_\calf\,s.
\end{equation}
\end{lemma}

By \eqref{Es-S-b}$_{1,2}$, the variations \eqref{E-sgeneral} preserve conformal foliations,
i.e., the property $\beta_{\mathcal D}\equiv0$.

\begin{lemma}[\textbf{Conservation laws}]\label{R-consumb}
Let $g_t\ (t\ge0)$ be ${\mathcal D}$-conformal
metrics \eqref{E-sgeneral} on a~foliated manifold $(M,\calf,{\mathcal D})$ such that
$H_0=-n\nabla^\calf\log u_0$ for a positive function $u_0\in C^\infty(M)$.
Then the following two functions and two vector fields on ${U}$ are $t$-independent\,$:$
\begin{equation*}
 \beta_{\mathcal D},\quad \|h_\calf\|^2/\|T\|,\quad
 H-(n/2)\nabla^\calf\log \|T\|,\quad
 H- n\nabla^\calf\log \|h_\calf\|.
\end{equation*}
\end{lemma}

\noindent\textbf{Proof}.
Using Lemma~\ref{L-btAt2} and $g^\perp(H,\,\cdot\,)=0$, we calculate
\begin{eqnarray*}
 \dt\,\|h\,\|^2\eq \dt\sum\nolimits_{\alpha}\tr(A_{{\mathcal E}_\alpha}^2)
 =2\sum\nolimits_\alpha\tr(A_{{\mathcal E}_\alpha}\dt A_{{\mathcal E}_\alpha})
 =-\sum\nolimits_\alpha {\mathcal E}_\alpha(s)\tr A_{{\mathcal E}_\alpha} = -g(\nabla s, H),\quad\\
 \dt g(H,H) \eq s\,g^\perp(H,H) +2\,g(\dt H, H) = -n\,g(\nabla s, H).
\end{eqnarray*}
Hence, $n\,\dt\beta_{\mathcal D}=\dt\,\|h\,\|^2_{g} -\frac1n\,\dt g(H,H)=0$,
that is the function $\beta_{\mathcal D}$ doesn't depend on $t$.

For any function $f\in C^1(M)$ and a vector $N\in {\mathcal D}_\calf$, using $(\dt g)(\cdot\,, N)=0$, we find
\[
 g(\nabla^\calf(\dt f), N) = N(\dt f)=\dt N(f)=\dt g(\nabla^\calf f, N) = g(\dt(\nabla^\calf f), N).
\]
Hence $\nabla^\calf(\dt f)=\dt(\nabla^\calf f)$.
By Lemma~\ref{L-btAt2}, we~find
\begin{equation}
\label{E-dtT}
 \dt\|T\|^2 =-\dt\sum\nolimits_{\alpha}\tr\big((T^\sharp_{{\mathcal E}_\alpha})^2\big)
 =-2\sum\nolimits_\alpha\tr(T^\sharp_{{\mathcal E}_\alpha}\dt T^\sharp_{{\mathcal E}_\alpha})
 =2s\sum\nolimits_\alpha\tr\big((T^\sharp_{{\mathcal E}_\alpha})^2\big) = -2\,s\,\|T\|^2.
\end{equation}
Similarly, by Lemma~\ref{R-consumb} and using the proof of Lemma~\ref{L-nablaNN}, we obtain
\begin{equation}\label{E-dthF}
 \dt\|h_\calf\|^2 = -s\,\|h_\calf\|^2.
\end{equation}
By the above, $h_\calf\ne0\ne T$ on ${U}$,
and we have $\dt\log\|T\|^2_{g_t} =-2\,s$ and $\dt\log\|h_\calf\|^2_{g_t} =-s$.
Using $\nabla^\calf\dt =\dt\nabla^\calf$, we obtain
 $\dt H_t =(n/2)\,\dt\nabla^\calf\log\|T\|_{g_t}$ and
 $\dt H_t =n\,\dt\nabla^\calf\log\|h_\calf\|_{g_t}$,
moreover, $\dt(\|h_\calf\|^2/\|T\|)=0$.
From the above the claim follows.
\qed

\subsection{Proofs of Propositions~\ref{T-main-loc}--\ref{T-mainA}, Theorem~\ref{T-main0} and Corollary~\ref{C-main1}}
\label{subsec:proofs}

\noindent\textbf{Proof of Proposition~\ref{T-main-loc}}.
Let $g_t=g_0+s\,g_0^\perp\ (0\le t<\eps)$ be
${\mathcal D}$-conformal
metrics on a~foliated manifold  $(M,\calf)$, where $s:M\times[0,\eps)\to\RR$ is a smooth function.
By Lemma~\ref{L-nablaNN}, $\calf$ is harmonic with respect to all $g_t$.
We differentiate (\ref{E-RicNs}) by $t$, and apply
 Lemmas~\ref{L-btAt2} and \ref{R-consumb} to obtain
\begin{equation*}
 \dt\Sc_{\,\rm mix}(g_t) = -(n/2)\,\Delta_\calf\,s +g(\nabla s, H)
 +s\big(\|h_\calf\|^2_{g_t} -2\,\|T\|^2_{g_t}\big).
\end{equation*}
Hence,
the linearization of (\ref{E-GF-Ricmix-mu}) at $g_0$
is the following linear PDE for $s$ on the leaves:
\begin{eqnarray*}
 \dt s =n\,\Delta_\calf\,s -2\,g_0(\nabla\,s,\,H_0)-2\,(\Sc_{\,\rm mix}(g_0)+\|h_\calf\|^2_{g_0}-2\,\|T\|^2_{g_0})s.
\end{eqnarray*}
The result follows from the theory of linear parabolic PDEs (see \cite{aub}) and the~finite holonomy assumption
(i.e., the local Reeb stability Theorem in Section~\ref{subsec:prel}).
\qed

\vskip1mm\noindent\textbf{Proof of Proposition~\ref{T-mainA}}.
By Theorem~\ref{T-main-loc}, (\ref{E-GF-Ricmix-mu}) admits a unit local leaf-wise smooth solution $g_t\ (0\le t<t_0)$.
The functions $\Sc_{\,\rm mix}(g_t), H_t$, $\|T\|_{g_t}$ and $\|h_\calf\|_{g_t}$ etc. are then uniquely determined for $0\le t<t_0$.
From \eqref{Es-S-b}$_3$ with $s=-2\,(\Sc_{\,\rm mix}(g) -\Phi)$ and using \eqref{E-RicNs} we obtain \eqref{Ec-ANtau1}.

($\bm i$)
By Lemma~\ref{R-consumb}(ii), $H_t-(n/4)\nabla^\calf\log\|T\|^2_{g_t}=X$
for some vector field $X$ on $M$. Since $H_0$ is conservative,
$X=-(n/4)\nabla^\calf\log\psi$ for some leaf-wise smooth function $\psi>0$ on $M$.
Hence,
 $H_t=-n\nabla^\calf\log\big(\psi^{1/4}\|T\|_{g_t}^{-1/2}\big)$
and, by condition $H_0=-n\nabla^\calf\log u_0$, one may take $\psi=u_0^4\,\|T\|^2_{g_0}$.
Define a leaf-wise smooth function $u:=(\Psi^\calf_2)^{1/4}\|T\|^{-1/2}_{g_t}$ on ${U}\times[0,t_0)$ and calculate
\[
 \dt(\log\|T\|^2_{g_t})=-4\,\dt\log(\|T\|_{g_t}^{-1/2})=-4\,\dt\log((\Psi^\calf_2)^{-1/4} u)=-4\,\dt\log u.
\]
By Lemma~\ref{R-consumb} and \eqref{E-Psi-i}
$\|h_\calf\|^2_{g_t}/\|T\|_{g_t}=\Psi^\calf_1(\Psi^\calf_2)^{-1/2}$,
thus, $u=(\Psi^\calf_1)^{1/2}\|h_\calf\|^{-1}_{g_t}$ on ${U}\times[0,t_0)$~and
\[
 \dt(\log\|h_\calf\|^2_{g_t})
 =-2\,\dt\log(\|h_\calf\|_{g_t}^{-1})=-2\,\dt\log((\Psi^\calf_1)^{-1/2} u)=-2\,\dt\log u.
\]
From the above and \eqref{E-RicNs} we then obtain
\begin{eqnarray*}
 \dt\log u\eq-(1/4)\,\dt(\log\|T\|^2_{g_t})=s/2=-\Sc_{\,\rm mix}(g_t)+\Phi\\
 \eq n\,\Delta_\calf\log u +n\,g(\nabla^\calf\log u, \nabla^\calf\log u) +n\,\beta_{\mathcal D}+\Phi
 +\Psi^\calf_1\,u^{-2}-\Psi^\calf_2\,u^{-4}.
\end{eqnarray*}
Substituting
 $\dt\log u={u}^{-1}\dt u$,
 $\nabla^\calf\log u={u}^{-1}\nabla^\calf u$ and
 $\Delta_\calf\log u={u}^{-1}\Delta_\calf\,u -{u}^{-2}g(\nabla^\calf u, \nabla^\calf u)$,
we find that $u$ solves the non-linear heat equation \eqref{Ec-dtrho1}.

($\bm{ii}$)
Note that $H$ obeys a forced leaf-wise Burgers equation {\rm (}a consequence of  \eqref{Ec-ANtau1}{\rm)}
\begin{equation*}
 \dt H +\nabla^\calf \|H\|_{g_t}^2= n\nabla^\calf(\Div_\calf H) -n^2\nabla^\calf\beta_{\mathcal D},
\end{equation*}
The rest of the proof see in \cite[Proposition~2]{rz}.
\qed

\vskip1mm\noindent\textbf{Proof of Theorem~\ref{T-main0}}.
(\textbf{i}) By Theorem~\ref{T-main-loc}, there exists a unique local solution $g_t$ on $M\times[0,t_0)$.
By Proposition~\ref{T-mainA}($ii$), $H$ obeys \eqref{Ec-ANtau1}, and $H=-n\nabla^\calf\log u$ for
some positive function $u$ satisfying \eqref{Ec-dtrho1} with $u(\cdot\,,0)=u_0$.
Note that conditions \eqref{E-cond1} yield $(u_0^-)^4\ge(\Psi^\calf_2)^+/(\Phi-n\lambda^\calf_0)$,
see \eqref{E-cond-u} with $\beta=\beta_{\mathcal D}+\Phi/n$ and $\lambda^\calf_0-\Phi/n<0$
and definitions  \eqref{E-Psi-i} and \eqref{dfmaxmin}.
By Theorem~\ref{prexistest},
one may leaf-wise smoothly extend a solution of \eqref{Ec-dtrho1} on $M\times[0,\infty)$,
hence $H_t(x)$ is defined for $t\ge0$ and is smooth on the leaves.
By Theorem~\ref{mainth1}($i$), $u\to\infty$ as $t\to\infty$ with exponential rate $n\alpha$ for
$\alpha\in(0,\min\{\lambda^\calf_1{-}\lambda^\calf_0,\,2(\Phi/n-\lambda^\calf_0\})$.
Hence, $\Psi^\calf_2 u^{-4}$ is leaf-wise smooth,
moreover, $\|T\|_{g_t}\to0$ and $h_\calf(g_t)\to0$ as $t\to\infty$.
 By~Theorem~\ref{mainth1}($ii$), $H_t=-n\nabla^\calf\log u$ approaches in $C^\infty$, as $t\to\infty$,
 to the vector field $\bar H=-n\nabla^\calf\log e_0$,
 hence $\Div_\calf\,H_t$ approaches to the leaf-wise smooth function $-n\,\Delta_\calf\log e_0$.
Since
 $-\Delta_\calf\,e_0-(\beta_{\mathcal D}+\Phi/n)\,e_0=\lambda^\calf_0\,e_0$,
we have, as $t\to\infty$,
\[
 \Div_\calf H_t-\|H\|^2_{g_t}/n\to
 -n(\Delta_\calf\,e_0)/e_0=n(\lambda^\calf_0+\beta_{\mathcal D})-\Phi.
\]
By \eqref{E-RicNs}, $\Sc_{\,\rm mix}(\cdot\,,t)$ approaches exponentially to $n\lambda^\calf_0-\Phi$ as $t\to\infty$.
Then a smooth solution to~(\ref{E-GF-Ricmix-mu}) is
$g_t=g_0\exp(-2\int_0^t(\Sc_{\,\rm mix}(\cdot\,,\tau)-\Phi)\dtau)$, where $t\ge0$,
see also Section~\ref{sec:cknorm}.

\vskip1mm
(\textbf{ii}) The smoothness of $g_t$ on $M$ follows from the finite holonomy assumption
and results of Section~\ref{sec:cknorm} (see also \cite{rz}).
Indeed, let $F$ be a leaf.
By the local Reeb stability Theorem (see Section~\ref{subsec:prel}), there is a~normal neighborhood
${\rm pr}:V\to F$ of $F$ (with a smooth normal section $Q$ -- an open $n$-dimensional disk) such that $(V,\calf_{\,|V},{\rm pr})$
is a foliated bundle diffeomorphic to $\widetilde F\times_{{\rm Hol}(F)} Q$.
There is a regular covering $\widetilde{\rm pr}:\tilde F\to F$ with covering group ${\rm Hol}(F)$,
since this group is finite, $\tilde F$ is a compact manifold.
The~normal neighborhood $V$, as a fiber bundle over $F$, can be pulled back
via $\widetilde{\rm pr}$ to a bundle $\tilde V$ over $\tilde F$ with the same fiber.
The standard pull-back construction yields a~canonical covering $\widetilde{\rm pr}: \tilde V\to V$ (a submersion),
enabling us to lift the foliation $\calf_{|V}$ to a~foliation $\tilde\calf$ of $\tilde V$,
transverse to the fibers and having $\tilde F$ as a leaf.
Let $\tilde x\in\tilde F$ with $\widetilde{\rm pr}(\tilde x) = x$.
Since the covering group is exactly the holonomy group of $F$,
$\widetilde{\rm pr}_* :\pi_1(\tilde F,\tilde x)\to\pi_1(F,x)$
injects $\pi_1(\tilde F,\tilde x)$ onto the subgroup of $\pi_1(F,x)$,
and the leaf $\tilde F$ of $\tilde\calf$ has trivial holonomy.
By~\cite[Theorem~2.4.1]{cc2}, $\tilde F$ has a neighborhood in $\tilde V$ that is a foliated disk bundle
with all leaves diffeomorphic to $\tilde F$.
 Let $\tilde g_t$ be the lifts of the metrics $g_t$
from $V$ onto the product $\tilde V$.
The~corresponding foliation on $\tilde V$ is harmonic,
the lift of the laplacian and potential function, $\tilde\beta_{\mathcal D}$, smoothly depend on $q$ on $\tilde V$.
(In case of totally geodesic foliation the leaves are isometric one to another,
see for example \cite{rov-m}.
Hence the leaf-wise laplacian on $\tilde V$ doesn't depend on~$q\in Q$).
The functions
$\tilde e_0$ and solution to \eqref{Ec-dtrho1}, $u_t\ (t>0)$, are smooth on $\tilde V$ (where $Q$ can be replaced by smaller normal section),
and they are lifts of leaf-wise smooth functions
$e_0$ and $u_t\ (t>0)$ on $V$ or on a~smaller neighborhood $W$ of $F$.
The vector fields $\tilde H_t$ and $\tilde{\bar H}$ are smooth on $\tilde V$,
and they are lifts of smooth vector fields $H_t$ and $\bar H$ on a neighborhood of $F$.
By the above, $g_t$ are smooth on~$M$.
\qed

\vskip1mm\noindent\textbf{Proof of Corollary~\ref{C-main1}}.
The claim ($i$) follows directly from Theorem~\ref{T-main0}.
The metrics $g_t\ (g_0=g)$ of Theorem~\ref{T-main0} diverge as $t\to\infty$ with the exponential rate
$\mu=\Phi-n\lambda^\calf_0$:
\[
 \exists\,C>1,\ \forall\,X\in {\mathcal D},\ \forall\,t\ge0:\
 C^{\,-1}e^{2\mu\,t}g(X,X)\le g_t(X,X)\le C e^{2\mu\,t}g(X,X).
\]
Consider ${\mathcal D}$-conformal metrics $\bar g_t=g_\calf + e^{-2\mu\,t} (g_t)^\perp$.
By \eqref{Es-S-b}$_3$, $\bar H_t=H_t$.
 Let $(\cdot\,,\,\cdot)_0$ be the inner product and the norm in $L_2(F)$ for any leaf $F$.
The function $v=e^{\,-\mu\,t}u$ converges as $t\to\infty$
to $\tilde u_0^0 e_0$, where $\tilde u_0^0=(\tilde u, e_0)_0=u_0^0+\int_0^\infty q_0(\tau)\dtau$,
see Theorem~\ref{mainth1} and \eqref{dftlu1}.
For $t\to\infty$ we have
\[
 \|h_\calf\|^2_{\bar g_t}=e^{2\mu\,t}\|h_\calf\|^2_{g_t}=\Psi^\calf_1/v^2\to\xi^2\|h_\calf\|^2_{g},\quad
 \|T\|^2_{\bar g_t}=e^{4\mu\,t}\|T\|^2_{g_t}=\Psi^\calf_2/v^4\to\xi^4\|T\|^2_{g},
\]
and the metrics $\bar g_t$ converge as $t\to\infty$ to the metric
$\bar g_\infty=g_\calf + \xi^{-2} g^\perp$.
By \eqref{E-RicNs}, we find
\[
 \Sc_{\,\rm mix}(\bar g_\infty)=n\lambda^\calf_0-\Phi +\xi^4\|T\|^2_{g} -\xi^2\|h_\calf\|^2_{g}\,.
\]
Comparing with \eqref{E-cond1} completes the proof of ($ii$).
\qed

\section{Results for PDEs}
\label{sec:Appendix}

The section plays an important role in proofs of main results (see Section~\ref{sec:egf}).

Let $(F, g)$ be a closed $p$-dimensional Riemannian manifold, e.g., a leaf of a foliation $\calf$.
Functional spaces over $F$ will be denoted without writing $(F)$, for example, $L_2$ instead of $L_2(F)$.

Let $H^l$ be the Hilbert space of differentiable by Sobolev real functions on $F$,
with the inner product $(\,\cdot,\cdot\,)_{l}$ and the norm $\|\cdot\|_l$.
In~particular, $H^0=L_2$ with the product $(\,\cdot,\cdot\,)_{0}$ and the norm~$\|\cdot\|_0$.

If $E$ is a Banach space, we denote by $\|\cdot\|_E$ the norm of vectors in this space.
If $B$ and $C$ are real Banach spaces, we denote by $\mathcal{B}^r(B,C)$ the Banach space
of all bounded $r$-linear operators $A:\,\prod_{i=1}^r B\rightarrow C$ with the norm
 $\|A\|_{\mathcal{B}^r(B,C)}=\sup_{v_1,\dots, v_r\in B\setminus{0}}\frac{\|A(v_1,\dots,v_r)\|_C}{\|v_1\|_B\cdot\ldots\cdot\|v_r\|_B}$.
If~$r=1$, we shall write $\mathcal{B}(B,C)$ and  $A(\cdot)$,
and if $B=C$ we shall write $\mathcal{B}^r(B)$ and $\mathcal{B}(B)$, respectively.

If $M$ is a $k$-regular manifold or an open neighborhood of the origin in a real Banach space, and $N$ is a real Banach space,
we denote by $C^k(M,N)$ $(k\ge 1)$ the Banach space of all $C^k$-regular functions $f:\,M\rightarrow N$,
for which the following norm is finite:
\begin{equation*}
 \|f\|_{C^k(M,N)}=\sup\nolimits_{x\in M}\max\{\|f(x)\|_N,\,\max\nolimits_{1\le |r|\le k}\|d^r f(x)\|_{\mathcal{B}^r(T_x M,\,N)}\}.
\end{equation*}
Denote by $\|\cdot\|_{C^k}$, where $0\le k<\infty$, the norm in the Banach space $C^k$; certainly, $\|\cdot\|_C$ when $k=0$.
In~coordinates $(x_1,\dots, x_p)$ on $F$, we have
$\|f\|_{C^k}=\max_{x\in F}\max_{|\alpha|\le k}|d^\alpha f(x)|$, where $\alpha\ge0$ is the multi-index of order
$|\alpha|=\sum_{i=1}^p \alpha_i$ and $d^\alpha$ is the partial derivative.

\subsection{The Schr\"{o}dinger operator}
\label{R-burgers-heat}

For a smooth
(non-constant in general) function $\beta$ on $(F,g)$, the Schr\"{o}dinger operator, see \eqref{E-Schr},
\begin{equation}\label{dfH}
 \mathcal{H}(u):=-\Delta u -\beta\,u
\end{equation}
is self-adjoint and bounded from below (but it is unbounded).
The domain of definition of $\mathcal{H}$ is~$H^2$.
 The spectrum $\sigma(\mathcal{H})$ of $\mathcal{H}$
 consists of an infinite sequence of isolated real eigenvalues
$\lambda_0\le\lambda_1\le\dots\lambda_j\le\dots$
counting their multiplicities, and $\lambda_j\rightarrow\infty$ as $j\rightarrow\infty$.
If we fix in $L_2$ an orthonormal basis of corresponding eigenfunctions $\{e_j\}$
(i.e., $\mathcal{H}(e_j)=\lambda_j e_j$) then any function $u\in L_2$ is expanded into the series
(converging to $u$ in the $L_2$-norm)
\begin{equation}\label{expan}
 u({x})=\sum\nolimits_{j=0}^\infty c_j\,e_j({x}),
 \qquad
 c_j = (u, e_j)_{0} =\int_{F} u(x)\,e_j(x)\dx.
\end{equation}
The proof of (\ref{expan}) is based on the following facts.
Since by the Elliptic regularity Theorem with $k=0$, we have $\mathcal{H}^{-1}: L_2\rightarrow H^2$ when
and the embedding of $H^2$ into $L_2$ is continuous and compact, see \cite{aub},
then the operator $\mathcal{H}^{-1}:\,L_2\rightarrow L_2$ is compact.
This means that the spectrum $\sigma(\mathcal{H})$
is discrete, hence by the spectral expansion theorem for compact self-adjoint operators,
$\{e_j\}_{j\ge0}$ form an orthonormal basis in $L_2$,
see \cite[I, Ch.~VII, Sect.~4; and II, Ch.~XII, Sect.~2]{ds}.
One can add a~constant to $\beta$ such that $\mathcal{H}$ becomes invertible in $L_2$
(e.g., $\lambda_0>0$)
and $\mathcal{H}^{-1}$ is bounded in $L_2$.

\begin{proposition}[see \cite{rz}]\label{P-lambda0-one}
Let $\beta$ be a smooth function on a closed Riemannian manifold $(F, g)$. Then the eigenspace of operator \eqref{dfH},
corresponding to the least eigenvalue, $\lambda_0$, is one-dimensional, and it contains a positive smooth eigenfunction, $e_0$.
\end{proposition}

The following facts will be used.

\vskip1mm\noindent
\textbf{Sobolev embedding Theorem} (see \cite{aub}).
\textit{If a nonnegative $k\in\ZZ$ and $l\in\NN$ are such that $2\,l>p+2\,k$,
then $H^l$ is continuously embedded into $C^k$.}

\vskip1mm\noindent
\textbf{Elliptic regularity Theorem} (see \cite{aub}).
\textit{If $\,\mathcal{H}$ is given by \eqref{dfH} and $0\notin\sigma(\mathcal{H})$,
then for any nonnegative $k\in\ZZ$ we have}
 $\mathcal{H}^{-1}:\; H^k\rightarrow H^{k+2}$.

\vskip1mm
The Cauchy's problem
for the \textit{heat equation with a linear reaction term}, see \eqref{Ec-dtrho1}, has a form
\begin{equation}\label{E-heat-in}
 \dt u = \Delta\,u +\beta\,u,\qquad u(x,0)=u_0(x).
\end{equation}
After scaling the time
and replacements of functions
\[
 t/n\to t,\quad
 \Psi^\calf_i/n\to\Psi_i,\quad
 \beta_{\mathcal D}+\Phi/n\to\beta,\quad
 \lambda^\calf_0-\Phi/n\to\lambda_0,
\]
the problem \eqref{Ec-dtrho1} reads as
the following Cauchy's problem for the \textit{non-linear heat equation} on~$(F,g)$:
\begin{equation}\label{Cauchy}
 \partial_t u=
 \Delta\,u +\beta\,u+\Psi_1(x)\,u^{-1}-\Psi_2(x)\,u^{-3},\qquad u(x,0)=u_0(x).
\end{equation}
By~\cite[Theorem~4.51]{aub}, the parabolic PDE \eqref{Cauchy} has a unique smooth solution
$u(\cdot\,,t)$ for $t\in[0,t_0)$.
Denote by $\mathcal{C}_{t}=F\times[0,t)$ the cylinder with the base $F$.
Define the quantities
\begin{eqnarray}\label{dfmaxmin}
\nonumber
 &&\Psi_i^+=\max\nolimits_{\,F}\,({\Psi_i}/{e_0^{2i}}),\quad
 \Psi_i^-=\min\nolimits_{\,F}\,({\Psi_i}/{e_0^{2i}}),\quad
 i=1,2,\\
 && u_0^+=\max\nolimits_{\,F}\,({u_0}/{e_0}),\quad
 u_0^-=\min\nolimits_{\,F}\,({u_0}/{e_0}),\quad
 \beta^-=\min\nolimits_{\,F}|\beta|.
\end{eqnarray}

We shall use the following scalar maximum principle \cite[Theorem~4.4]{ck1}.

\begin{proposition}\label{P-weak-max}
Suppose that $X(t)$ is a smooth family of
vector field on a closed Riemannian manifold $(F,g)$,
and $f\in C^\infty(\RR\times[0, T))$. Let $u : F\times[0, T)\to\RR$ be a $C^\infty$ supersolution to
\[
 \dt u \ge\Delta_{g} u +\<X(t),\nabla u\> +f(u,t).
\]
Suppose that $\varphi:[0,T]\to\RR$ solves the Cauchy's problem for ODEs $(d/dt)\varphi = f(\phi(t),t),\ \varphi_1(0)=C$.
If $u(\cdot, 0)\ge C$ then $u(\cdot, t)\ge\varphi(t)$ for all $t\in[0, T)$.
(Claim also holds with the sense of all three inequalities reversed).
\end{proposition}

\subsection{Long-time solution}
\label{sec:existest}


\begin{lemma}\label{L-apriori}
Let $\lambda_0<0$ for $(F,g)$ and $u(x,t)>0$ be a solution in $\mathcal{C}_{t_0}$ of
\eqref{Cauchy} with the condition
\begin{equation}\label{E-cond-u}
 (u_0^-)^4 \ge {\Psi^+_2}/{|\lambda_0|},
\end{equation}
see \eqref{E-cond1}. Then the following a priori estimates are valid:
\begin{equation}\label{estCauchy}
 w_{-}(t)\le u(x,t)/e_0(x)\le w_{+}(t),\qquad (x,t)\in\mathcal{C}_{t_0},
\end{equation}
where $\lambda_0<0$ and
\begin{equation}\label{dfwmM}
 w_{-}(t)=e^{-\lambda_0 t}\big((u_{0}^-)^4+\frac{\Psi^+_2}{\lambda_0}
 -e^{\,4\lambda_0 t}\frac{\Psi^+_2}{\lambda_0}\big)^{\frac{1}{4}},\quad
 w_{+}(t)= e^{-\lambda_0 t}\big(
  (u_{0}^+)^2 -\frac{\Psi^+_1}{\lambda_0} +e^{\,2\lambda_0 t}\frac{\Psi^+_1}{\lambda_0}
 \big)^{\frac{1}{2}}.
\end{equation}
\end{lemma}

\noindent\textbf{Proof}. Since $e_0(x)>0$ on $F$, we can change the unknown function in (\ref{Cauchy}):
\begin{equation*}
 u(x,t)=e_0(x)\,w(x,t).
\end{equation*}
Substituting into (\ref{Cauchy})
and using $\Delta e_0+\beta e_0=-\lambda_0e_0$, we obtain the Cauchy's problem for $w(x,t)$:
\begin{equation}\label{Cauchw}
\partial_t w=\Delta w-\lambda_0 w+2\,g(\nabla\log e_0,\,\nabla w)
 +e_0^{-2}(x){\Psi_1(x)}\,{w^{-1}} -e_0^{-4}(x){\Psi_2(x)}\,{w^{-3}},\quad
 w(\cdot\,,0)={u_0}/{e_0}.
\end{equation}
Then, using \eqref{dfmaxmin}$_{1}$, we obtain the differential inequalities
\begin{equation*}
 \Delta w-\lambda_0 w+2\,g(\nabla\log e_0,\nabla w)-{\Psi^+_2}{w^{-3}}\le\partial_t w
 \le\Delta w -\lambda_0 w+2\,g(\nabla\log e_0,\nabla w)+{\Psi^+_1}{w^{-1}}\,.
\end{equation*}
By the scalar maximum principle of Proposition~\ref{P-weak-max} and
\eqref{dfmaxmin}$_{2,3}$, we conclude that \eqref{estCauchy} holds,
 where
$w_{-}(t)$ and $w_{+}(t)$ are solutions of the following Cauchy's  problems for ODEs
\begin{equation*}
 \frac{\rm d}{\rm dt}\,w_{-}=-\lambda_0\,w_{-}-{\Psi^+_2}w_{-}^{-3},\quad w_{-}(0)=u_0^-;\qquad
  \frac{\rm d}{\rm dt}\,w_{+}=-\lambda_0\,w_{+}+{\Psi^+_1}w_{+}^{-1},\quad w_{+}(0)=u_0^+.
\end{equation*}
One may check that these solutions are expressed by (\ref{dfwmM})
and $w_{-}(t)<w_{+}(t)$ for all $t\ge0$.
\qed

\vskip1mm
Note that if $\Psi^+_i=0$ (i.e., $\Psi_i\equiv0$) then \eqref{estCauchy} reads
$u_{0}^-e^{-\lambda_0 t} \le u(\cdot,t)/e_0\le u_{0}^+e^{-\lambda_0 t}$.
Define
\[
 v(x,t)=e^{\lambda_0 t}u(x,t),
\]
see \eqref{Cauchy}, and obtain the Cauchy's problem
\begin{equation}\label{Cauchv1}
 \partial_t v= \Delta\,v +(\beta+\lambda_0)\,v + Q,
 \quad v(x,0)=u_0(x),
\end{equation}
where
 $Q:=\sum_{i=1}^2(-1)^{i+1}\Psi_i(x)\,v^{1-2i}(x,t)\,e^{\,2i\lambda_0 t}$.
Certainly, $Q=Q_1-Q_2$, where
\begin{equation*}
 Q_1(x,t) = \Psi_1(x)v^{-1}(x,t)e^{2\lambda_0 t},\quad
 Q_2(x,t) = \Psi_2(x)v^{-3}(x,t)e^{4\lambda_0 t}.
\end{equation*}

\begin{lemma}\label{lmestder1}
Let $v(x,t)$ be a positive solution of \eqref{Cauchv1} in $\mathcal{C}_{t_0}=F\times[0,\,{t_0})$,
where $\lambda_0<0$, the functions $u_0>0$ and $\Psi_i\ge0$ belong to $C^\infty$ and \eqref{E-cond-u} is satisfied.
Then

$($i$)$~for any multi-index $\alpha=(\alpha_1,\alpha_2,\dots,\alpha_p)$ there
exists a real $C_\alpha\ge0$ such that
\[
 |\partial_x^\alpha v(x,t)|\le C_\alpha(1+t)^{|\alpha|},\quad
 (x,t)\in\mathcal{C}_{t_0}.
\]

$($ii$)$~for any multi-index $\alpha$ there exist real $\bar Q_{i\alpha}\ge0\ (i=1,2)$ such that
\begin{equation}\label{estderQ}
 |\partial_x^\alpha Q_i(x,t)|\le\bar Q_{i\alpha}(1+t)^{|\alpha|}e^{\,2\,i\lambda_0 t},
 \quad
 (x,t)\in\mathcal{C}_{t_0},\quad i=1,2.
\end{equation}
\end{lemma}

\noindent\textbf{Proof}.
Using \eqref{estCauchy} and \eqref{dfwmM}, we estimate the solution $v(x,t)$
of \eqref{Cauchv1} when \eqref{E-cond-u} holds:
\begin{equation}\label{estv1}
 v_-\le{v(x,t)}/{e_0(x)}\le v_+,\qquad (x,t)\in F\times[0,\infty),
\end{equation}
where the constants are given by
 $v_-=\big((u_0^-)^4-{\Psi^+_2}/{|\lambda_0|}\big)^{\frac{1}{4}}$ and
 $v_+=\big((u_0^+)^2+{\Psi^+_1}/{|\lambda_0|}\big)^{\frac{1}{2}}$.

($\bm{i}$) Denote for brevity by $D_j=\partial_{x_j}\,(j=1,2,\dots,p)$ and
$D_\alpha=\partial_x^\alpha=D_{\alpha_1}D_{\alpha_2}\cdots D_{\alpha_p}$.
Differentiating \eqref{Cauchv1} by $x_1,\dots,x_p$,
we obtain the following PDEs for the functions $p_\alpha(x,t):=\partial_x^\alpha v(x,t)$:
\begin{eqnarray}\label{eqpj}
\nonumber
 \partial_t p_j\eq(\Delta +(\lambda_0+\beta)\id)p_j +D_j(\beta) v +D_j(Q), \\
 \partial_t p_{jk}\eq(\Delta {+}(\lambda_0{+}\beta)\id)p_{jk}
 +D_j(\beta)p_k +D_k(\beta)p_j +D_{jk}(\beta) v +D_{jk}(Q),
\end{eqnarray}
and so on, where $1\le j,k\le p$ and
\begin{eqnarray}\label{dfaxt}
\nonumber
 D_j(Q) \eq \sum\nolimits_{i} (-1)^{i+1}e^{\,2i\lambda_0t} v^{-2i}\big(D_j(\Psi_i)v +(1-2\,i)\Psi_i p_j\big),\\
\nonumber
 D_{jk}(Q) \eq \sum\nolimits_{i} (-1)^{i+1}e^{\,2i\lambda_0t} v^{-2i}\big(
 D_{jk}(\Psi_i)v +(1-2\,i)D_{j}(\Psi_i) p_k +(1-2\,i)D_{k}(\Psi_i) p_j\\
 \minus 2\,i(1-2\,i) v^{-1}\Psi_i p_j p_k +(1-2\,i)\Psi_i p_{jk}\big),
\end{eqnarray}
and so on. Let us change  unknown functions in equations (\ref{eqpj}),
and so on:
\begin{equation}\label{chngvarpjk}
 p_j=\tilde p_j\,e_0,\quad p_{jk}=\tilde p_{jk}\,e_0,\ \ldots
\end{equation}
Then in the same manner, as (\ref{Cauchw}) have been
obtained from (\ref{Cauchy}), we get for $j,k=1,2,\dots,p$
\begin{eqnarray}\label{eqtlpj}
\nonumber
 \partial_t\tilde p_j\eq\Delta\tilde p_j+2\,g(\nabla\log e_0,\nabla\tilde p_j)
 +a\tilde p_j+{b_j}/{e_0} +D_j(\beta ) v/e_0,\\
 \partial_t\tilde p_{jk}\eq\Delta\tilde p_{jk} +2\,g(\nabla\log e_0,\nabla\tilde p_{jk})
 +a\tilde p_{jk}+{b_{jk}}/{e_0} +D_{jk}(\beta) v/e_0,
\end{eqnarray}
and so on, where
\begin{eqnarray*}
 a \eq\sum\nolimits_i (-1)^{i+1}(1-2\,i)\Psi_i\,v^{-2i}e^{2i\lambda_0t},\quad
 b_j=\sum\nolimits_i (-1)^{i+1} D_j(\Psi_i)\,v^{1-2i}e^{\,2i\lambda_0t},\\
 b_{jk} \eq \sum\nolimits_i (-1)^{i+1}e^{\,2i\lambda_0t}v^{-2i}
 \big( D_{jk}(\Psi_i)\frac v{e_0} +(1-2\,i)\big(
 D_j(\Psi_i)\tilde p_k +D_k(\Psi_i)\tilde p_j -2\,i\,\Psi_i\frac{e_0}v\,\tilde p_j\,\tilde p_k \big) \big).
\end{eqnarray*}
From (\ref{estv1}) and (\ref{dfaxt})\,--\,(\ref{eqtlpj}) we get the differential inequalities
\begin{eqnarray*}
 &&
 -a^+(t)|\tilde p_j| -b_j^+ -\beta_j^+ v_+ \le
 \partial_t\tilde p_j -\Delta\tilde p_j -2\,g(\nabla\log e_0,\nabla\tilde p_j)
 \le a^+(t)|\tilde p_j|+ b_j^+ +\beta_j^+ v_+\,,\\
 &&-a^+(t)|\tilde p_{jk}| -b_{jk}^+ -\beta_{jk}^+ v_+ \le
 \partial_t\tilde p_{jk} -\Delta\tilde p_{jk} -2\,g(\nabla\log e_0,\nabla\tilde p_{jk})
 \le a^+(t)|\tilde p_{jk}|+ b_{jk}^+ +\beta_{jk}^+ v_+
\end{eqnarray*}
for $j=1,2,\dots,p$, where
\begin{eqnarray}\label{dfaplt}
\nonumber
 &&\hskip-10mm
 a^+(t)=\sum\nolimits_i (2\,i-1)\Psi^+_i (v_-)^{-2i}e^{2i\lambda_0t},\quad
 b_j^+=\sum\nolimits_i \big((v_-)^{1-2i}\max_F \big|D_j(\Psi_i)\big|\big),\quad
 \beta_j^+=\max_{\,F}\big|D_j(\beta)\big|,\\
\nonumber
 &&\hskip-10mm
 b_{jk}^+=
 \sum\nolimits_i e^{\,2i\lambda_0t}(v_-)^{-2i}
 \big(\max_F \big|D_{jk}(\Psi_i)/e_0^{2i}\big|\,v_+
 +(2\,i-1)\big(
 \max_F \big|D_j(\Psi_i)/e_0^{2i}\big|\,\tilde p_k \\
 && +\max_F \big|D_k(\Psi_i)/e_0^{2i}\big|\,\tilde p_j
 +2\,i\,(v_-)^{-1}\Psi_i^+\tilde p_j\,\tilde p_k \big) \big),\quad
 \beta_{jk}^+ = \max_{\,F}\big|D_{jk}(\beta)\big|.
\end{eqnarray}
By the maximum principle of Proposition~\ref{P-weak-max}, the estimate $|\tilde p_j(x,t)|\le\tilde p_j^+(t)$
is valid for any $(x,t)\in\mathcal{C}_\infty=F\times[0,\infty)$,
where $p_j^+(t)$ solves the Cauchy's problem for the ODE:
\begin{equation*}
 \frac{\rm d}{\rm dt}\,p_j^+=a^+(t)| p_j^+| + b_j^+ +\beta_j^+ v_+,
 \quad p_j^+(0)=\bar p_j^0 := \max\nolimits_{\,F}\big|\tilde p_j(\cdot\,,0)\big|\,.
\end{equation*}
 As is known,
\begin{equation*}
 p_j^+(t)=\bar p_j^0\exp\big(\int_0^ta^+(\tau)\dtau\big)
 +\int_0^t (b_j^+ +\beta_j^+ v_+) \exp\big(\int_s^ta^+(\tau)\dtau\big)
 \ds\,.
\end{equation*}
In view of (\ref{dfaplt})$_1$, we have
\begin{equation*}
 \int_0^\infty a^+(\tau)\dtau<\infty,
\end{equation*}
the above yield that for any $j\in\{1,2,\dots,p\}$ there exists a real $\tilde C_j>0$ such that
\begin{equation*}
 |\tilde p_j(x,t)|\le\tilde C_j(1+t),\qquad
 (x,t)\in\mathcal{C}_\infty\,.
\end{equation*}
In view of (\ref{chngvarpjk}), this completes the proof of (${i}$) for $|\alpha|=1$.

Similarly we obtain that for any
$j,k\in\{1,2,\dots,p\}$ there exists a real $\tilde C_{jk}\ge0$ such that
 $|\tilde p_{jk}(x,t)|\le\tilde C_{jk}(1+t)^2$
for $(x,t)\in\mathcal{C}_\infty$.
By (\ref{chngvarpjk}), we obtain claim (${i}$) for $|\alpha|=2$.
By induction with respect to $|\alpha|$ we prove (${i}$) for any~$\alpha$.

($\bm{ii}$) Estimates (\ref{estderQ}) for $|\alpha|=0$,
 $|Q_i(x,t)|\le (\max\nolimits_{\,F}\Psi_i)\,(v_-)^{1-2i}e^{\,2i\lambda_0 t}$,
follow immediately from (\ref{estv1}).
Estimates (\ref{estderQ}) for $|\alpha|=1,2$ follow from claim (${i}$),
estimates (\ref{estv1}) and equalities \eqref{dfaxt}.
By~induction with respect to $|\alpha|$ we prove (${ii}$) for any~$\alpha$.
\qed

\begin{theorem}\label{prexistest}
The Cauchy's problem \eqref{Cauchy} on $F$, with $\lambda_0<0$ and the
initial value $u_0(x)$ satisfying \eqref{E-cond-u}, admits a~unique
smooth solution $u(x,t)>0$ in the cylinder $\mathcal{C}_\infty=F\times[0,\infty)$.
\end{theorem}

\noindent\textbf{Proof}.
The positive solution $u(x,t)$ of  \eqref{Cauchy} satisfies a priory estimates (\ref{estCauchy})
on any cylinder $\mathcal{C}_{t_\star}$ where it exists.
By standard arguments, using the local theorem of the existence and uniqueness for semi-flows,
we obtain that this solution can be uniquely prolonged on the cylinder $\mathcal{C}_\infty$.
Then, by Lemma~\ref{lmestder1}, all partial derivatives by $x$ of $u(x,t)$ exist in $\mathcal{C}_\infty$.
Hence, $u$ is smooth on $\mathcal{C}_\infty$.
\qed

\subsection{Asymptotic behavior of solutions}
\label{sec:asympttestshr}

Recall that $\lambda_0$ and $e_0>0$ are the least eigenvalue and the ground state of
the operator \eqref{dfH}.

\begin{theorem}\label{mainth1}
Let $u>0$ be a smooth solution on $\mathcal{C}_\infty$ of
\eqref{Cauchy} with $\lambda_0<0$ and the initial value $u_0(x)$
satisfying \eqref{E-cond-u} $($see Theorem~$\ref{prexistest})$.
Then there exists a
solution $\tilde u$ on
$\mathcal{C}_\infty$ of the linear PDE
\begin{equation}\label{Schrheat}
 \partial_t\tilde u=\Delta\tilde u+(\beta(x)+\lambda_0)\,\tilde u
\end{equation}
such that for any $\alpha\in\big(0,\min\{\lambda_1-\lambda_0,\,2\,|\lambda_0|\}\big)$
and any $k\in\NN$
\[
({i}) \ u=e^{-\lambda_0 t}(\tilde u+\theta(x,t)), \qquad
({ii}) \ \nabla\log u=\nabla\log e_0+\theta_1(x,t),
\]
where $\|\theta(\cdot\,,t)\|_{C^k}=O(e^{-\alpha t})$
and
$\|\theta_1(\cdot\,,t)\|_{C^k}=O(e^{-\alpha t})$ as $t\rightarrow\infty$.
\end{theorem}

\noindent\textbf{Proof}.
($\bm{i}$) Let $G_0(t, x, y)$ be the fundamental solution of (\ref{Schrheat}), called the \textit{heat kernel}.
As is known,
$G_0(t, x, y)=\sum\nolimits_{j} e^{(\lambda_0-\lambda_j)\,t}e_j(x)\,e_j(y)$.
Due to the Duhamel's principle, the solution $v=e^{\lambda_0 t}u$ of the Cauchy's
problem (\ref{Cauchv1}) satisfies the nonlinear integral equation
\begin{equation}\label{duhamel1}
 v(x,t)=\int_{F} G_0(t,x,y)\,u_0(y)\dy
 +\int_0^t\Big(\int_{F} G_0(t-\tau, x,y)Q(y,\tau)\dy\Big)\dtau\,.
\end{equation}
Expand $v$, $u_0$ and $Q$
into Fourier series by eigensystem~$\{e_j\}$:
\begin{equation}\label{expv1}
\begin{array}{cc}
  v(x,t) =\sum\nolimits_{j=0}^\infty v_j(t)\,e_j(x), & v_j(t)= (v(\cdot\,,t),e_j)_0= \int_{F}v(y,t)\,e_j(y)\dy, \\
  u_0(x) =\sum\nolimits_{j=0}^\infty u_j^0\,e_j(x), & u_j^0= (u_0,e_j)_0= \int_{F}u_0(y)\,e_j(y)\dy, \\
  Q(x,t) =\sum\nolimits_{j=0}^\infty q_j(t)\,e_j(x), & q_j(t)= (Q(\cdot\,,t),e_j)_0= \int_{F}Q(y,t)\,e_j(y)\dy.
\end{array}
\end{equation}
Then we obtain from (\ref{duhamel1}):
\begin{equation}\label{frmvj1}
 v_j(t)=u_j^0\,e^{(\lambda_0-\lambda_j)t}+\int_0^t e^{(\lambda_0-\lambda_j)(t-\tau)}q_j(\tau)\dtau
 \qquad (j=0,1,\dots).
\end{equation}
Substituting $v_j(t)$ of \eqref{frmvj1} into (\ref{expv1}), we represent $v$ in the form
 $v=\tilde u+\theta$, where
\begin{eqnarray}\label{dftlu1}
 \tilde u\eq \tilde u_0^0\,e_0 +\sum\nolimits_{j=1}^\infty u_j^0\,e^{(\lambda_0-\lambda_j)t}e_j,\quad
 \tilde u_0^0=u_0^0+\int_0^\infty q_0(\tau)\dtau,\\
\label{dftht1}
 \theta\eq-\Big(\int_t^\infty q_0(\tau)\dtau\Big) e_0+\sum\nolimits_{j=1}^\infty\tilde v_j\,e_j,\qquad
 \tilde v_j=\int_0^te^{(\lambda_0-\lambda_j)(t-\tau)}q_j(\tau)\dtau.
\end{eqnarray}
Observe that $\tilde u$ solves (\ref{Schrheat}) with the initial condition
 $\tilde u(\,\cdot\,,\,0)=u_0+\big(\int_0^\infty q_0(\tau)\dtau\big)e_0$.

Let us take $k\in\NN$, $l=\big[p/4+k/2\big]+1$ and
$\gamma<\lambda_0$. Using assumption $u_0\in C^\infty(F)$ and the
fact that $Q(\cdot\,,t)\in C^\infty(F)$ for any $t\ge 0$, we may consider the functions
 $w_0:=(\mathcal{H}-\gamma\id)^lu_0$ and
 $P(\,\cdot\,,t):=(\mathcal{H}-\gamma\id)^lQ(\,\cdot\,,t)$,
which have the same properties: $w_0\in C^\infty(F)$ and $P(\cdot\,,t)\in C^\infty(F)$ for any $t\ge 0$.
Let us represent
\begin{eqnarray*}
 (u_0,\,e_j)_0\,e_j\eq((\mathcal{H}-\gamma\id)^{-l}w_0,\,e_j)_0\,e_j
 =(w_0,\,(\mathcal{H}-\gamma\id)^{-l}e_j)_0\,e_j\\
 \eq(w_0,\,e_j)_0\frac{e_j}{\lambda_j-\gamma}=(\mathcal{H}-\gamma\id)^{-l}\big((w_0,\,e_j)_0\,e_j\big).
\end{eqnarray*}
Similarly, we obtain
\[
 (Q(\,\cdot\,,t),e_j)_0\,e_j=(\mathcal{H}-\gamma\id)^{-l}\big((P(\,\cdot\,,t),e_j)_0\,e_j\big).
\]
Using (\ref{dftht1}) and taking into account that the operator
$(\mathcal{H}-\gamma\id)^{-l}$ acts continuously in $L_2$ and
that the series in (\ref{dftlu1}) and (\ref{dftht1}) converge in
$L_2$, we obtain the representations:
\begin{eqnarray*}
 \sum\nolimits_{j=1}^\infty u_j^0e^{(\lambda_0-\lambda_j) t}e_j \eq(\mathcal{H}-\gamma\id)^{-l}
 \sum\nolimits_{j=1}^\infty e^{(\lambda_0-\lambda_j) t}(w_0,e_j)_0\,e_j,\\
 \sum\nolimits_{j=1}^\infty\tilde v_j(t)e_j\eq (\mathcal{H}-\gamma\id)^{-l}
 \int_0^t\big(\sum\nolimits_{j=1}^\infty e^{(\lambda_0-\lambda_j)(t-\tau)}(P(\,\cdot\,,t),e_j)_0\,e_j\big)\dtau.
\end{eqnarray*}
By the Elliptic Regularity Theorem and
the Sobolev Embedding Theorem (see Section~\ref{R-burgers-heat}), the operator $(\mathcal{H}-\gamma\id)^{-l}$
acts continuously from $L_2$ into $C^k$. Then we have
\begin{eqnarray}\label{estsum1}
 &&\hskip-8mm\|\sum\nolimits_{j=1}^\infty u_j^0e^{(\lambda_0-\lambda_j) t}e_j\|_{C^k}
 \le\|(\mathcal{H}-\gamma\id)^{-l}\|_{\mathcal{B}(L_2,C^k)}\cdot
 \|\sum\nolimits_{j=1}^\infty e^{(\lambda_0-\lambda_j) t}(w_0,e_j)_0\,e_j\|_0\\
\nonumber
 &&\hskip-4mm =\|(\mathcal{H}{-}\gamma\id)^{-l}\|_{\mathcal{B}(L_2,C^k)}
 \big(\sum\nolimits_{j=1}^\infty e^{2(\lambda_0-\lambda_j) t}(w_0,e_j)_0^2\big)^{1/2}
 \le\|(\mathcal{H}{-}\gamma\id)^{-l}\|_{\mathcal{B}(L_2,C^k)}\,e^{(\lambda_0-\lambda_1) t}\|w_0\|_0,\\
\label{estsum2}
 &&\hskip-8mm \|\sum\nolimits_{j=1}^\infty\tilde v_j(t)e_j\|_{C^k}
 \le\|(\mathcal{H}-\gamma\id)^{-l}\|_{\mathcal{B}(L_2,C^k)}\cdot\|\int_0^t\sum\nolimits_{j=1}^\infty
 e^{(\lambda_0-\lambda_j)(t-\tau)}(P(\,\cdot\,,t),e_j)_0\,e_j\dtau\|_0\nonumber\\
 &&\hskip-4mm\le\|(\mathcal{H}-\gamma\id)^{-l}\|_{\mathcal{B}(L_2,C^k)}\int_0^t\|\sum\nolimits_{j=1}^\infty
 e^{(\lambda_0-\lambda_j)(t-\tau)}(P(\,\cdot\,,t),e_j)_0\,e_j\|_0\dtau\nonumber\\
 &&\hskip-4mm\le\|(\mathcal{H}-\gamma\id)^{-l}\|_{\mathcal{B}(L_2,C^k)}\int_0^t
 \big(\sum\nolimits_{j=1}^\infty e^{2(\lambda_0-\lambda_j)(t-\tau)}(P(\,\cdot\,,t),e_j)_0^2\big)^{1/2}
 \dtau\nonumber\\
 &&\hskip-4mm\le\|(\mathcal{H}-\gamma\id)^{-l}\|_{\mathcal{B}(L_2,C^k)}
 \int_0^te^{(\lambda_0-\lambda_1)(t-\tau)}\|P(\,\cdot\,,t)\|_0\dtau.
\end{eqnarray}
On the other hand, by Lemma~\ref{lmestder1}(ii),
\begin{equation*}
 \|P(\,\cdot\,,t)\|_0\le\sqrt{{\rm Vol}(F,g)}\,\|(\mathcal{H}-\gamma\id)^l\big(Q(\,\cdot\,,t)\big)\|_{C^0}
 \le\bar Q(1+t)^{2l}e^{2\lambda_0t}
\end{equation*}
for some $\bar Q\ge0$. Then, continuing \eqref{estsum2}, we find
\begin{eqnarray}\label{estintP}
 && \int_0^te^{(\lambda_0-\lambda_1)(t-\tau)}\|P(\,\cdot\,,t)\|_0\dtau
 \le\bar Q\int_0^te^{(\lambda_0-\lambda_1)(t-\tau)}(1+\tau)^{2l}e^{2\lambda_0\tau}\dtau\\
\nonumber
 &&\hskip-8mm <\bar Q e^{(\lambda_0-\lambda_1) t}(1+t)^{2l}\int_0^t e^{(\lambda_1-\lambda_0+2\lambda_0)\,\tau}\dtau
 =\bar Q(1+t)^{2l}
 \Big\{\begin{array}{ll}
 \frac{e^{2\lambda_0 t}
 -e^{(\lambda_0-\lambda_1)t}}{\lambda_1-\lambda_0+2\,\lambda_0}&\rm{if}\quad 2\lambda_0\ne
 \lambda_0-\lambda_1,\\
 e^{(\lambda_0-\lambda_1)t} t
 \phantom{\frac{\frac iI}I}
 &\rm{if}\quad 2\lambda_0=\lambda_0-\lambda_1.
 \end{array}
\end{eqnarray}
From (\ref{dftlu1})\,--\,(\ref{estintP}) we get claim (${i}$).

\noindent\ ($\bm{ii}$)
From (\ref{dftlu1}) and (\ref{dftht1}) we obtain
\begin{equation*}
 u=e^{-\lambda_0 t}\big(\tilde u_0^0\,e_0+\bar\theta(\cdot\,,t)\big),\quad
 \nabla u=e^{-\lambda_0 t}\big(\tilde u_0^0\nabla e_0+\nabla\bar\theta(\cdot\,,t)\big),
\end{equation*}
where $\bar\theta(\cdot\,,t)=\theta(\cdot\,,t)+\sum\nolimits_{j=1}^\infty u_j^0e^{(\lambda_0-\lambda_j) t}e_j$.
In view of (\ref{estsum1}), $\|\bar\theta(\cdot\,,t)\|_{C^k}=O(e^{-\alpha t})$
for any $k\in\NN$.
Furthermore, since $\tilde u(\cdot,0)>0$ on $F$, then $\tilde
u_0^0=(\tilde u(\,\cdot,\,0),\;e_0)_0>0$.
Using
\[
 w(\cdot,t,\tau):=\tau u(\cdot,t)+(1-\tau)\,\tilde u_0^0\,e^{-\lambda_0t}e_0
 =e^{-\lambda_0 t}(\tilde u_0^0\,e_0+\tau\,\bar\theta(\cdot,t)),
\]
we have
\begin{eqnarray*}
 \theta_1(\cdot,t)\eq\nabla\log u(\cdot,t)-\nabla\log e_0
 =\int_0^1\frac{\partial}{\partial\tau}\Big(\nabla\log w(\cdot,t,\tau)\Big)\dtau\\
 \eq\int_0^1\Big(\frac{\nabla\bar\theta(\cdot,t)}{\tilde u_0^0 e_0+\tau\,\bar\theta(\cdot,t)}
 -\frac{\bar\theta(\cdot,t)(\tilde u_0^0\nabla e_0+\tau\nabla\bar\theta(\cdot,t))}
 {(\tilde u_0^0 e_0+\tau\,\bar\theta(\cdot,t))^2}\Big)\dtau\,.
\end{eqnarray*}
By the above, and the fact that
 $\inf\{\,|\tilde u_0^0 e_0+\tau\,\bar\theta(\cdot,t)|:\ x\in F,\ t\in[t_0,\infty),\ \tau\in[0,1]\}>0$
holds for $t_0>0$ large enough, follows claim (${ii}$).
\qed

\subsection{Nonlinear heat equation with parameter}
\label{sec:cknorm}

Let the metric $g$, the connection $\nabla$ and the laplacian $\Delta$ smoothly depend on $q$, which
belongs to an open subset $Q$ of $\RR^n$.
Consider the Cauchy's problem
on a closed Riemannian manifold $(F,g)$:
\begin{equation}\label{nnlinheat}
 \partial_tu=\Delta u+f(x,u,q),\quad u(x,0,q)=u_0(x,q).
\end{equation}
Here, $f$ is defined in the domain
$D=F\times I\times Q$, where $I\subseteq\RR$ is an
interval, and $u_0$ is defined in the domain $\tilde D=F\times
Q$ and satisfies the condition: $u_0(x,q)\in I$ for any $x\in F$ and $q\in Q$.

\begin{proposition}
Suppose that $f\in C^\infty(D)$, $u_0\in C^\infty(\tilde D)$,
all partial derivatives of $f$ and $u_0$ by $x$, $u$ and $q$ are
bounded in $D$ and $\tilde D$, and for any $q\in Q$ there exists an
unique solution $u:F\times[0,T]\times Q\to\RR$ of the Cauchy's problem \eqref{nnlinheat}
such that all its partial derivatives by $x$ are bounded in $F\times[0,T]\times Q$.
Then $u(\cdot\,,t,\cdot)\in C^\infty(F\times Q)$ for any $t\in[0,T]$.
\end{proposition}

\noindent\textbf{Proof} is standard, we give it for the convenience of a reader.
As is known, $u(\cdot\,,t,q)\in C^\infty(F)$ for any $q\in Q$, $t\in[0,T]$.
We should prove the smooth dependence on $q$ of the
solution $u(x,t,q)$ and of all its partial derivatives by $x$ for
any fixed $t\in[0,T]$. We shall divide the proof into several steps.

\textbf{Step 1}:
the continuous dependence of $u(x,t,q)$ in $q$.
To show this, take $q_0\in Q$ and denote by
$v(x,t,q)=u(x,t,q)-u(x,t,q_0)$ and $v_0(x,q)=u_0(x,q)-u_0(x,q_0)$. Let
us represent:
\begin{equation*}
 f(x,u(x,t,q),q)-f(x,u(x,t,q_0),q_0)=F(x,t,q)v(x,t,q)+G(x,t,q)\cdot(q-q_0),
\end{equation*}
where
\begin{equation*}
\begin{array}{c}
 F(x,t,q)=\int_0^1\partial_uf(x,\,u(x,t,q_0)+\tau v(x,t,q),\,q_0+\tau(q-q_0))\dtau \\
 \qquad G(x,t,q)=\int_0^1{\rm grad}_qf(x,\,u(x,t,q_0)+\tau v(x,t,q),\,q_0+\tau(q-q_0))\dtau\,.
\end{array}
\end{equation*}
Then the function $v(x,t,q)$ is a solution of the Cauchy's problem:
\begin{equation*}
\partial_t v=\Delta v+F(x,t,q)v+G(x,t,q)\cdot(q-q_0),\quad
v_{\,\vert\,t=0}=G_0(x,q)\cdot(q-q_0),
\end{equation*}
where
 $G_0(x,q) =\int_0^1{\rm grad}_q v_0(x,\,q_0+\tau(q-q_0))\dtau$.
Then by the maximum principle of Proposition~\ref{P-weak-max},
\begin{equation}\label{maxprinc}
 |v(x,t,q)|\le w(t,q)\qquad
 \forall\,(x,t,q)\in F\times[0,T]\times Q,
\end{equation}
where $w(t,q)$ is the solution of the following Cauchy's problem for the ODE:
\begin{equation}\label{ODE}
 \dt\,w=\bar F|w|+\bar G\,|q-q_0|,\qquad w(0,q)=\bar G_0\,|q-q_0|
\end{equation}
with
 $\bar F=\sup\limits_{F\times[0,T]\times Q}|F|$,
 $\bar G=\sup\limits_{F\times[0,T]\times Q}|G|$
 and
 $\bar G_0=\sup\limits_{F\times Q}|G_0|$.
Then, from (\ref{maxprinc}) and (\ref{ODE}) we get
\begin{equation*}
 |v(x,t,q)|\le\big(\bar G_0 e^{\bar Ft}+(e^{\bar Ft}-1)({\bar G}/{\bar F})\big)|q-q_0|,
 \qquad(x,t,q)\in F\times[0,T]\times Q,
\end{equation*}
which implies the claim of Step 1.

\textbf{Step 2}: all the partial derivatives of $u(x,t,q)$
by $x$ are continuous in $q$. Differentiating subsequently by
$x$ both sides of the equation and of the initial condition in
(\ref{nnlinheat}), we have the following Cauchy's problems for
$p_\alpha=\partial^\alpha_xu$ ($\alpha$ is the multi-index):
\begin{eqnarray}\label{Cauchpi}
\nonumber
 \partial_t p_i\eq\Delta p_i +\partial_u
 f(x,u(x,t,q),q)p_i+\partial_x^if(x,u(x,t,q),q),\quad
 p_{i\,\vert\,t=0}=\partial_x^iu_0(x,q),\\
\nonumber
 \partial_t p_{i,j}\eq\Delta p_{i,j} +\partial_u f(x,u(x,t,q),q)p_{i,j}+\partial^2_u
 f(x,u(x,t,q),q)p_ip_j\\
 \plus\partial_u\partial_x^if(x,u(x,t,q),q)p_j+\partial_x^{i,j}f(x,u(x,t,q),q),\quad
 p_{i,j\,\vert\,t=0}=\partial_x^{i,j}u_0(x,q),
\end{eqnarray}
and so on. Applying the claim of Step 1 to these Cauchy's problems, we
prove the claim of Step 2.

\textbf{Step 3}:
$u(x,t,q)$ is smooth with respect to $q$.
Take $q_0\in Q$ and consider the divided difference
\[
 \delta_{\,{\bf y}}(x,t,s)=\frac1{s}\,(u(x,t,q_0+s{\bf y})-u(x,t,q_0))\quad
 ({\bf y}\in\RR^n,\;s\in\RR).
\]
Denote by $\delta^0_{\,{\bf y}}(x,s)=\frac{u_0(x,q_0+s{\bf y})-u_0(x,q_0)}{s}$.
As in the Step 1, we obtain the Cauchy's problem for $\delta_{\,{\bf y}}$
\begin{eqnarray}\label{Cauchy4}
 \partial_t\delta_{\,{\bf y}}=\Delta \delta_{\,{\bf y}}+\tilde F(x,t,s)\delta_{\,{\bf y}}
 +\tilde G(x,t,s)\cdot{\bf y},
 \quad \delta_{\,{\bf y}\,\vert\,t=0}=\delta^0_{\,{\bf y}}(x,s),\\
 \nonumber
\begin{array}{c}
 \tilde F(x,t,s) =\int_0^1\partial_uf(x,u(x,t,q_0)+\tau(u(x,t,q_0+s{\bf y})
 -u(x,t,q_0))),\ q_0+s{\bf y})\dtau,\\
 \tilde G(x,t,s) =\int_0^1{\rm grad}_q f(x,u(x,t,q_0)+\tau(u(x,t,q_0+s{\bf y})-u(x,t,q_0))),\ q_0+s{\bf y})\dtau.
\end{array}
\end{eqnarray}
Applying to the Cauchy's problem (\ref{Cauchy4}) the claim of Step 1,
we conclude that $\delta_{\,{\bf y}}(x,t,s)$ is continuous by $s$ at the
point $s=0$, that is there exists the directional derivative
 $d_{\,{\bf y}}(x,t,q_0)={\rm grad}_q u(x,t,q_0)\cdot{\bf y}=\lim\nolimits_{\,s\rightarrow 0}\delta_{\,{\bf y}}(x,t,s)$.
Moreover, $d_{\,{\bf y}}(x,t,q)$ is the solution of Cauchy's problem
\begin{eqnarray}\label{Cauchy5}
\nonumber
 \partial_t d_{\,{\bf y}}\eq\Delta d_{\,{\bf y}}+\partial_u f(x,u(x,t,q),q)d_{\,{\bf y}}
 +{\rm grad}_qf(x,u(x,t,q)\cdot {\bf y},\\
 d_{\,{\bf y}\,\vert\,t=0}\eq{\rm grad}_q u_0(x,q)\cdot{\bf y}.
\end{eqnarray}
Applying the claim of Step 1 to this Cauchy's problem, we find
that $d_{\,{\bf y}}(x,t,q)={\rm grad}_q u(x,t,q)\cdot{\bf y}$
continuously depends on $q$ for any ${\bf y}\in\RR^n$.
Thus, $u(x,t,q)$ is $C^1$-regular in $q$.
Applying the above arguments to the Cauchy's problem (\ref{Cauchy5}),
we conclude that $u(x,t,q)$ belongs to $C^2$ with respect to~$q$.
Finally, we prove by induction that $u(x,t,q)$ is smooth in~$q$.

\textbf{Step 4}. Applying all the arguments of Step 3 to the Cauchy's problems
(\ref{Cauchpi}) and so on, we prove that all derivatives
of $u(x,t,q)$ in $x$ smoothly depend on $q$.
\qed

\end{document}